\newcommand{\beq}{\begin{equation}}   
\newcommand{\eeq}{\end{equation}}

\documentclass[a4paper]{article}

\usepackage{indentfirst}
\usepackage{graphicx,amssymb} %

\date{} %

\usepackage{lipsum}

\newcommand\blfootnote[1]{%
  \begingroup
  \renewcommand\thefootnote{}\footnote{#1}%
  \addtocounter{footnote}{-1}%
  \endgroup
}

\usepackage{amsmath}
\usepackage[utf8]{inputenc}
\usepackage[english]{babel}
\usepackage{amsthm}

\newtheorem{theorem}{Theorem}[section]

\newtheorem{proposition}[theorem]{Proposition}
\newtheorem{corollary}[theorem]{Corollary}
\newtheorem{lemma}[theorem]{Lemma}

\newtheorem*{rmk}{Remark}
\usepackage{amssymb}
\usepackage{amsfonts}
\usepackage{tikz-cd}

\begin{document}
\title{Brill--Noether theory and Green's conjecture for general curves on simple abelian surfaces}

\author{Federico Moretti}

\maketitle

\thispagestyle{empty}


\begin{abstract}
We compute the gonality and the dimension of the Brill--Noether loci $W^1_d(C)$ for curves in a non-primitive linear system of a simple abelian surface, adapting vector bundles techniques \`a la Lazarsfeld originally introduced with $K3$ surfaces. As a corollary, we obtain Green's conjecture for general curves on abelian surfaces.

\end{abstract}

\begin{center}
  {\textit{Dedicated to my Grandparents}}
  \end{center}

\section{Introduction}\blfootnote{F. MORETTI, \\ Humboldt Universität zu Berlin, Institüt für Mathematik, Unter den Linden 6, 10099 Berlin. \emph{Current adress:} State University of New York at Stony Brook, 100 Niccols Road, 11794 (NY), USA. \\
e-mail: federico.moretti@stonybrook.edu}In this paper, we compute the dimension of the Brill--Noether loci $\mathcal{W}^1_d(|L|)$ for non-primitive linear system on abelian surfaces. We establish a linear growth condition for them and, as a corollary, we obtain Green's conjecture for general curves in $|L|$.

Brill--Noether theory has been deeply investigated for curves lying on $K3$ surfaces. Lazarsfeld proved that a  smooth curve in the primitive linear system of a Picard general $K3$ surface is Brill--Noether general \cite{LAZ}, giving a proof of the Brill--Noether--Petri conjecture. Green and Lazarsfeld proved that the Clifford index of a smooth curve on a $K3$ is constant when the curve varies in its linear system. Moreover, when the Clifford index is not maximal it is computed by the restriction of a line bundle defined on the surface; see \cite{GL87}. Later,  Ciliberto, Pareschi and Knutsen proved that the gonality is always constant in all such linear systems, with some well-understood exceptions (see \cite{CILI},\cite{Kn09}). For the invariance of the Clifford index and for the fact that $\mathcal{O}_S(C)_{|_C}\simeq \omega_C$, curves on $K3$ surfaces turned out to be natural candidates to satisfy Green's conjecture on the syzygies of canonical curves
\[
K_{p,2}(C,\omega_C)=0 \quad \mathrm{for} \quad p<\mathrm{Cliff}(C),
\]
where $K_{i,j}(C,\omega_C)$ denotes the $(i,j)$-th Koszul cohomology group of the canonical bundle $\omega_C$. Voisin proved Green's conjecture for an arbitrary smooth curve in the primitive linear system of a general $K3$ (see \cite{V1},\cite{V2}). Aprodu and Farkas later generalized Voisin's work, proving Green's conjecture for all curves lying on any $K3$ surface in \cite{FA} (see also \cite{FApp},\cite{FAP},\cite{KEM2} for recent proofs of general Green's conjecture with new insights). We should point out here that this circle of ideas found other applications in the study of Brill--Noether theory and Green's conjecture for curves lying on certain surfaces with non positive canonical bundle. See for instance \cite{KL1},\cite{KL2} for Enriques surfaces and \cite{L} for Del-Pezzo surfaces (and other rational surfaces). Notice that all these surfaces satisfies $h^1(\mathcal{O}_S)=0$.

Abelian surfaces share some common behavior with $K3$ surfaces. It is natural to ask to what extent a similar study for their curves sections can be carried out. In recent years Knutsen, Lelli-Chiesa and Mongardi started studying Brill--Noether theory for curves in the primitive linear system $|H|$ of a general abelian surface \cite{KLL}. Among other things, they obtained that the general curve in $|H|$ is Brill--Noether general, but (unlike in the $K3$ case) the Clifford index is not constant for smooth curves in $|H|$. In particular, there are curves carrying linear series with negative Brill--Noether number and they can prove under certain hypotheses that these Brill--Noether loci are of the expected codimension in $|H|$.  A similar result was obtained before with different techniques by Paris in  \cite{PAR} (unpublished) and generalized later by Bayer and Li in \cite{BAYER}.

In this paper, we start the study of Brill--Noether theory for curves on a non-primitive linear system on a simple abelian surface (i.e. there is no elliptic curve $E \subset S$) and its implications for Green's conjecture. We present here our results for abelian surfaces whose N\'eron-Severi is of rank $1$.
Here are the main results of this paper.
\begin{theorem}
  Let $(S,N)$ be a $(1,e)$ polarized abelian surface with $\mathrm{NS}(S)=\mathbb{Z}N$ and $e\ge 2$. Let $C$ be a general curve of genus $g=m^2e+1$ in $|mN|$. If $m \le 2$ then $C$ is of maximal gonality. Otherwise the following hold:
  \begin{eqnarray}
    \mathrm{gon}(C)&=&2e(m-1)+2. \nonumber \\
    \mathrm{dim}W^1_d(C)&\le& d-2e(m-1)-2 \quad \quad \mathrm{for} \quad d\le g-2e(m-1). \nonumber 
    \end{eqnarray}
    Moreover, the general curve carries at least $\frac{N^2(N^2-1)}{2}$ minimal pencils.
    \label{t1}
\end{theorem}
The case $m=1$ follows from \cite{KLL},\cite{BAYER},\cite{PAR} and it is not tractable with the methods presented in the paper.
In \cite[Theorem 2]{Ap} Aprodu established a sufficient condition for Green's conjecture in terms of the dimension of the Brill--Noether loci $W^1_d(C)$; this yields the following:

\begin{theorem}
  Let $(S,N)$ be a simple polarized abelian surface as above. Then the general element of $|mN|$ satisfies Green's conjecture for any $m>0$.
  \label{t2}
\end{theorem}

In section 3, while exhibiting a minimal pencil on the general curve in $|mN|$, we will observe a non-general Brill--Noether behaviour (for $m>2$) in $\mathcal{G}^1_{\mathrm{2e(m-1)+2}}$. Indeed, we will prove that for general $C$ in $|L|$ the Brill--Noether locus $W^1_{\mathrm{gon}(C)}(C)$ is non-reduced. This fact makes the existence of the minimal pencil non-trivial and poses additional challenges with respect to the case of $K3$ surfaces (see section $3$). This is the most serious difference with the $K3$ case.

Moreover, the components of lower gonality of $|mN|$ have a component of higher codimension than expected. We remark that the above is different from what happens in the primitive linear system.

The curves on non-primitive linear systems on an abelian surface together with a torsion line bundle coming from the surface are special even from the Prym point of view. Indeed, from the previous Theorem, the following Corollary follows immediately.
\begin{corollary}
  Let $(C,\eta_{|_{C}})\in \mathcal{R}_g$ be a general curve in $|mN|$ of gonality $2e(m-1)+2$ together with a torsion line bundle $\eta \in \mathrm{Pic}^0S$ with $\eta^{\otimes 2}=\mathcal{O}_S$. Then the \'etale cover $\tilde{C} \to C$ induced by $\eta$ is of gonality $\mathrm{gon}(\tilde{C})\le 4e(m-1)+2< 2\mathrm{gon}(C)$. Indeed, there is a commutative diagram
  \[
  \begin{tikzcd}
    \tilde{C} \arrow{r} \arrow{d} & C \arrow{d} \\
    \tilde{S} \arrow{r} & S
    \end{tikzcd}
  \]
  where $\tilde{S} \to S$ is the $(1,2e)$ abelian surface induced by  $\eta$.
  \end{corollary}

In the following, we sum up some differences and similarities with linear systems on $K3$ surfaces.
\begin{corollary}
  Let $(S,N)$ as above and $L=mN$. The following hold:
  \begin{itemize}
  \item The Clifford index is not constant among the curves in $|L|$;
  \item The Clifford index is never computed by the restriction of a line bundle on $S$;
  \item For $m\ge 3$ the map computing the gonality $C \to \mathbb{P}^1$ can be chosen to be the restriction of a rational map $S \dashrightarrow \mathbb{P}^1$ for the general $C \in |L|$. More precisely, the Clifford index and the gonality can be computed by a line subbundle of $N\otimes \eta _{|_C}$ for some $\eta \in \mathrm{Pic}^0S$.
    \end{itemize}
  \end{corollary}

  Given the above corollary, there is a natural connection with a conjecture due to Donagi and Morrison (see  \cite[conjecture 1.2]{DM}) for curves on $K3$ surfaces: for a curve $C$ lying on an abelian surface do all the special $g^r_d$ with $d\le g-1$ arise as subbundles of the restriction of a line bundle on the surface of degree $\le g-1$? The above corollary gives some evidence for that.
  Donagi-Morrison's conjecture holds for $K3$ surfaces under suitable conditions, see \cite{LCDM}.  Studying the conjecture for pencils on curves on abelian surfaces may be an interesting problem. 
  
  Before giving the outline of the paper, we introduce a more general setting and some notation. $S$ will now denote a simple abelian surface and $C \subset S$ a smooth curve of genus $g$, with $L=\mathcal{O}_S(C)$. Moreover we define
\[
g(S,L)=\mathrm{min}\big{\{}\{M\cdot (L-M) +2 | \quad M \in \mathrm{Pic}(S), 2 \le h^0(L-M) \le h^0(M)\} \cup \{ \lfloor \frac{g+3}{2} \rfloor \}\big{\}}.
\]
Abusing notation $N(S,L) \in \mathrm{Pic}(S)$ will denote a divisor achieving the above minimum (i.e. $N(S,L)\cdot( L-N(S,L))=g(S,L)-2$ and $2 \le h^0(N(S,L))\le h^0(L \otimes N(S,L)^\vee)$). If $\mathrm{NS}(S)=\mathbb Z \cdot N$ and $L=mN$ we get $g(S,L)=2e(m-1)+2$ and $N(S,L)=N$ as long as $m\ge 2$ and $N^2\ge 4$.

Let us state a more general version of Theorem \ref{t2}. When writing $M> N$ we mean that the divisor $M-N$ is effective ($M-N\ge 0$ includes the possibility of $M-N$ being numerically trivial). Moreover, we will often abuse notation denoting with the same symbol a divisor and the associated line bundle.
\begin{theorem}
  Let $S$ be a simple abelian surface and  $|L|$ a complete linear system such that  $N(S,L)$ exists. Suppose moreover that either $L\ge 3N(S,L)$ or $\mathrm{NS}(S)=\mathbb Z \cdot N(S,L)$. Then, for any $d\le g-g(S,L)+2$ and any irreducible component $\mathcal{W} \subset \mathcal{W}^1_d|L|$ dominating $|L|$, the following inequality holds:
  \[
  \mathrm{dim}(\mathcal{W}) \le d-g(S,L)+g-2.
  \]
  Moreover, for any  smooth $C \in |L|$ it holds
  \[
  g(S,L)-2\le \mathrm{gon}(C) \le g(S,L)
  \]
  and $g(S,L)$ is the gonality of the general curve in the linear system. Furthermore, the locus of curves of gonality $g(S,L)-i$ have a component of codimension $i$, for $i=1,2$.
  \label{t4}
\end{theorem}

Let us remark that in case  $g(S,L)=  \frac{g+3}{2} $ the linear growth condition still makes sense since (in this situation) it is required up to degree $g(S,L)-1$.  The case of the primitive linear system of a general abelian surface $(S,L)$ despite being included in Theorem \ref{t1} is not included in Theorem \ref{t4} since it cannot be deduced with the methods of the present paper; it follows directly from the more general fact that in this case the general curve in $|L|$ is Brill--Noether general (see \cite{PAR}, \cite{KLL}, \cite{BAYER}).
Notice that the first condition in the above Theorem is satisfied in all but finitely many linear systems on a simple abelian surface. As before, we obtain as a corollary that the general curve in a linear system meeting one of the extra conditions listed in the Theorem satisfies Green's conjecture.
 We remark that the hypothesis that the abelian surface is simple is unavoidable for Theorem \ref{t1} as the following example shows. Consider $S =E \times E$ where $E$ is an elliptic curve, the general curve $C \in |2E\times\{0\}+\{0\} \times nE|$ is smooth bielliptic of genus $2n+1$ so cannot satisfy Theorem \ref{t1} (it has infinitely many minimal pencils). If we consider sufficiently ample linear system, the upper bound for the gonality still works, but it is far from being sharp. In order to prove Theorem \ref{t2} one should follow other paths (for instance Green's conjecture is known for the general cover of an elliptic curve, see \cite{KEM}).

  The easiest cases we left out of consideration are the linear systems of type $|L|=|H_1 + nH_2|$ in abelian surfaces such that $\mathbb{Z}H_1 \oplus \mathbb{Z}H_2\subset \mathrm{NS}(S)$ (with $H_1,H_2$ primitive) and $n$ is small. In this case, we cannot prove that the general element $C \in |L|$ carries at least one $g^1_{g(S,L)}$, but we have evidence that the bounds on the dimension of $\mathcal{W}^1_d|L|$ still work. It would be interesting to study this case, in principle, there may be new non trivial examples of curves carrying a one-dimensional family of minimal pencils.

  {\bf{Plan of the paper.}} The paper is organised as follows: in section 2 we give an upper bound for the dimension of the space of linear series of curves in a given linear system in the surface $\mathcal{W}^1_d(|L|)$. In section 3 we give an upper bound for the gonality obtained by specializing line bundles coming from the surface, this will prove the first bound to be sharp and will give good control on the growth of the dimension of the Brill--Noether loci.
 One studies pencils on the general curve $ C\in |L|$ via the associated Lazarsfeld--Mukai bundle (see the beginning of the following section). Because one can always realize such bundles as an extension of positive rank $1$ torsion-free sheaf, we are able to control the dimension of these extensions and a fortiori to give an upper bound of the dimension of $\mathcal{W}^1_d|L|$. The argument to show that such bundles fit in such an extension follows from a deformation theoretic argument of Pareschi when $h^1(E)=0$ (see Section 2.1.1) and it requires some extra vector bundles arguments when $h^1(E)>0$ (see section 2.1.2). 
In order to compute the gonality, it is then enough to exhibit a pencil of degree $g(S,L)$. We construct a scheme $F$ parametrizing triples $(\varphi,(P_1,\dots,P_{N(S,L)^2-2}),C)$, where $\varphi:S \dashrightarrow \mathbb{P}^1 $ is a rational map induced by $2$ sections of $H^0(N(S,L)\otimes \eta)$ (for $\eta \in\mathrm{Pic}^0(S)$), $(P_1,\dots,P_{N^(S,L)^2-2})$ are base points of the rational map and $C \in H^0(L \otimes \mathcal{I}_{P_1,\dots,P_{N(S,L)^2-2}})$. We have a canonical morphism $F \to \mathcal{W}^1_{g(S,L)}|L|$ given by
\[
(\varphi,(P_1,\dots,P_{N(S,L)^2-2}),C) \to (C, \varphi_{|_C}).
\]
Our result is complete once we show that the composition \[F \to \mathcal{W}^1_{g(S,L)}|L| \to |L|\] is dominant.

Ideas in section $2$ mainly come from \cite{FA} (with technical adjustments since $h^1(\mathcal{O}_S)=2$, see for instance lemma \ref{diff}, section 2.1.2, or the proof of \ref{oie}), section $3$ consists of direct computations in Brill--Noether theory and section 4 is a study of certain subschemes of $\mathrm{Hilb}^l(S)$ arising from Lazarsfeld--Mukai bundles. Let us remark here that for the lower bound the differences from \cite{FA} are mainly technical and can all be avoided if one is just interested in the lower bound for the gonality (and not in the full linear growth condition). The most serious difference with the $K3$ case is the existence of a minimal pencil of the expected degree for the general curve in the linear system; the fact that the tangent space to $W^1_{g(S,L)}(C)$ at a pencil coming from the surface is always non-trivial poses additional challenges (one needs to rule out the possibility that the pencils coming from the surfaces may induce a $1-$ dimensional $W^1_{g(S,L)}(C)$ for $C$ varying in a codimension $1$ family).

{\bf{Acknowledgements.}} I would like to thank my advisor Gavril Farkas for introducing me to the topic and for many useful discussions as well as for his support. I would also like to thank Andrei Bud for his careful reading and his comments on the first draft. Finally, I am thankful to Margherita Lelli-Chiesa for reading the paper and for some useful suggestions.
I also benefited from many useful discussions with Andrea Di Lorenzo.
This project has received funding from the European Research Council (ERC) under the
European Union Horizon 2020 research and innovation program (grant agreement No.
834172).

The author thanks the referee for careful reports that remarkably improved the exposition.

\section{The dimension of $\mathcal{W}^1_d(L)$}

From now on $\mathcal{W}\subset\mathcal{W}^1_d|L|$ will denote an irreducible component dominating $|L|$. We recall the definition of $N(S,L)\in \{ N \in \mathrm{Pic}(S) | \quad 2 \le h^0(N)\le h^0(L-N)\}$ from the introduction, i.e. a divisor (not unique) minimizing the following quantity:
\[
N\cdot (L-N).
\]
We will focus on linear systems such that $N(S,L)$ exists. By definition \[g(S,L)=N(S,L)\cdot (L - N(S,L))+2.\] 
The goal of this section is to prove the following result:
\begin{theorem}
  Let $S$ be a simple abelian surface and $|L|$ a non-trivial system on $S$ such that $L\ge 3 N(S,L)$ or $L=2N(S,L)$ and $\mathrm{NS}(S)=\mathbb Z\cdot N(S,L)$ , then a dominating component $\mathcal{W} \subset \mathcal{W}^{1}_d|L|$  satisfies the following linear growth condition: for  $d \le g-g(S,L)+2$ we have
  \[
  \mathrm{dim}\mathcal{W}\le d-g(S,L)+g-2.
  \]
  \label{t3}
\end{theorem}
Notice in particular that in case $L=2N(S,L)$ and $\mathrm{NS}(S)=\mathbb Z \cdot N(S,L)$ the above already implies that the general curve in $|L|$ is of maximal gonality.
For technical reasons, we will consider the algebraic system $|L|_{num}$, including all the translates $|L \otimes \eta|$ ($\eta \in \mathrm{Pic}^0(S)$).
The thesis of the theorem is equivalent to $\mathrm{dim}\mathcal{W}\le d-g(S,L)+g$  for any dominating irreducible component  $\mathcal{W}\subset \mathrm{dim}\mathcal{W}^1_d(|L|_{num})$.  

We recall the construction of Lazarsfeld--Mukai bundles $E_{C,A}$ associated to any base point free pencil $(C,A) \in {W}^1_d|C|$. The vector bundle $E^\vee_{C,A}$ is defined as the kernel of the evaluation map $\mathrm{ev}$
\[
\begin{tikzcd}
0 \arrow{r} & E_{C,A}^\vee \arrow{r} & H^0(A) \otimes \mathcal{O}_S \arrow{r}{ev} & A \arrow{r} & 0.
\end{tikzcd}
\]
Dually
\[
\begin{tikzcd}
  0 \arrow{r} & H^0(A)^\vee \otimes \mathcal{O}_S \arrow{r} & E_{C,A} \arrow{r} & \omega_C \otimes A^\vee \arrow{r} & 0.
  \end{tikzcd}
\]
One gets $c_1(E_{C,A})=L,c_2(E_{C,A})=d$.
We have the following easy.
\begin{lemma}
\label{0}
  Let $S,C,A$ be as above then
  \begin{eqnarray}
    0 &\le& h^1(E_{C,A}) \le 4; \nonumber \\
    h^0(E_{C,A})&=&g-1-d+h^1(E_{C,A}) \ge 3. \nonumber
  \end{eqnarray}

  \end{lemma}
\begin{proof}
    For the last inequality notice that $g-1-d\ge g(S,L)-3\ge N(S,L)^2-1\ge 3 $.
\end{proof}

\subsection{Structure of Lazarsfeld--Mukai bundles}
We are interested in an upper bound of the dimension of an irreducible component of $\mathcal{W}^1_d|L|$ dominating $|L|$.
In this section, we show that we may suppose that for a general element $(C,A)$ in such a component, the associated Lazarsfeld--Mukai bundle fits in a short exact sequence \[
\begin{tikzcd}
  0 \arrow{r} & M \arrow{r} & E_{C,A} \arrow{r} & N \otimes \mathcal{I}_\xi \arrow{r} & 0,
  \end{tikzcd}
\]
with $h^0(M)\ge 2$. We will sometimes write $E=E_{C,A}$ for brevity.

  \subsubsection{Case $h^1(E)=0$.}

This is well known for $K3$ surfaces (see for instance \cite{pare}). We argue on the same lines of \cite{FA}.  Let $A \in W^r_d(C)-W^{r+1}_d(C)$ be a globally generated line bundle, we consider the Petri map
    \[
    \begin{tikzcd}
      \mu_{0,A}:H^0(C,A) \otimes H^0(C,\omega_C \otimes A^\vee) \arrow{r} & H^0(C,\omega_C),
      \end{tikzcd}
    \]
    whose kernel can be described in terms of Lazarsfeld--Mukai bundles. Indeed, we define $M_A$ as
    the rank $r$ vector bundle on $C$ defined as the kernel of the evaluation map
    \[
    \begin{tikzcd}
      0 \arrow{r} & M_A \arrow{r} & H^0(A) \otimes \mathcal{O}_C \arrow{r}{ev} & A \arrow{r} & 0.
      \end{tikzcd}
    \]
    Twisting the above with $\omega_C \otimes A^\vee$ we get $\mathrm{ker}(\mu_{0,A})=H^0(C, M_A \otimes \omega_C \otimes A^\vee)$. There is also an exact sequence on $C$:
    \[
    \begin{tikzcd}
      0 \arrow{r} & \mathcal{O}_C \arrow{r} & E_{C,A}^\vee \otimes \omega_C \otimes A^\vee \arrow{r} & M_A \otimes \omega_C \otimes A^\vee \arrow{r} & 0.
    \end{tikzcd}
    \]
    On the other hand twisting the defining sequence of $E_{C,A}$ with $E_{C,A}^\vee$ we get
    \[
    \begin{tikzcd}
      0 \arrow{r} & H^0(A)^\vee \otimes E_{C,A}^\vee \arrow{r} & E_{C,A} \otimes E_{C,A}^\vee \arrow{r} & \omega_{C} \otimes A^\vee \otimes E_{C,A}^\vee \arrow{r} & 0.
      \end{tikzcd}
    \]
    Now if $h^1(E)=0$ the map $H^0(E_{C,A} \otimes E_{C,A}^\vee )\to H^0( \omega_{C} \otimes A^\vee \otimes E_{C,A}^\vee)$ is surjective. We deduce
    \[
    H^0(E_{C,A} \otimes E_{C,A}^\vee )= H^0(\omega_{C} \otimes A^\vee \otimes E_{C,A}^\vee).
    \]
    Now we recall a lemma from \cite{pare}, which follows from Sard's lemma applied to $\mathcal{W}^r_d|L| \to |L|$.
    \begin{lemma}
      Suppose $\mathcal{W} \subset \mathcal{W}^r_d|L|$ is a dominating component, and $(C,A) \in \mathcal{W}$ is a general element such that $A$ is globally generated and $h^0(A)=r+1$. Then the coboundary map \\ $H^0(C, M_A \otimes \omega_C \otimes A^\vee) \to H^1(C ,\mathcal{O}_C)$ is zero.
      \label{sbattiiii}
      \end{lemma}
    So we get (see also \cite{ApP08}, Corollary 3.3):
    \begin{proposition}
      If $\mathcal{W} \subset \mathcal{W}^1_d|L|$ is a dominating component, and $(C,A) \in \mathcal{W}$ is a general element such that $A$ is globally generated, $h^0(A)=2$, $h^1(E)=0$; then \\ $\mathrm{dim}W^1_d(C) \le \rho(g,1,d)+h^0(C,E_{C,A} \otimes E_{C,A}^\vee)-1$. Moreover, equality holds if and only if $\mathcal{W}$ is reduced at $(C,A)$.

      \label{313}
    \end{proposition}
    \begin{proof}
        The orthogonal of the tangent space of $W^1_d(C)$ at a  $A\in W^1_d(C)$ can be canonically 
        identified with the image of the Petri map. Hence
        \[
        \mathrm{dim}(T_AW^1_d(A))= g-(r+1)(g-d+r)+h^0(M_A\otimes \omega_C\otimes A^\vee)=\rho(g,r,d)+h^0(M_A\otimes \omega_C\otimes A^\vee).
        \]
        If $A$ is globally generated we get, by Lemma \ref{sbattiiii}, \[
       h^0(M_A\otimes \omega_C\otimes A^\vee)=h^0(E_{C,A}\otimes E_{C,A}^\vee)-1,
        \]
        the Proposition follows.
    \end{proof}
    Hence either the general $E_{C,A}$ has an endomorphism $\varphi$ realizing $E_{C,A}$ as an element of \\ $ \mathrm{Ext}^1(\mathrm{Im}(\varphi),\mathrm{ker}(\varphi))$ (one may suppose that $\varphi$ drops rank everywhere) or the dimension of $W^1_d(C)$ is $\rho(g,1,d)$. Indeed up to replacing a non trivial $\varphi$ with $\varphi-\lambda\cdot \mathrm{id}$ we may suppose that $\varphi$ drops rank at a point and hence everywhere (the determinant is a section of $\mathcal O_S$). Now both $\mathrm{ker}(\varphi)$ and $\mathrm{Im}(\varphi)$ inject into $E$ and at least one of them must have $2$ sections (since $h^0(E)\ge 3)$, taking the double dual we obtain the desired line subbundle with $2$ linearly independent global sections.

    \subsubsection{Case $h^1(E)>0$.}
   As in the previous section the task is to find a subline bundle of $E$ with two linearly independent global sections. As a consequence $E$ will automatically lie in a sequence of the form\[
\begin{tikzcd}
  0 \arrow{r} & M \arrow{r} & E_{C,A} \arrow{r} & N \otimes \mathcal{I}_\xi \arrow{r} & 0,
  \end{tikzcd}
\] 
up to replacing the subline bundle with two global sections with a more positive one in order to remove the torsion at the cokernel.

Since $h^1(E)>0$ we may consider a nontrivial extension of vector bundles
   \[
  \begin{tikzcd}
    0 \arrow{r} & \mathcal{O}_S \arrow{r} & F \arrow{r} & E \arrow{r} & 0.
    \end{tikzcd}
  \]
  We fix $M=N(S,L)$ \footnote{This choice of notation because we are looking for the subline bundle $M\to E$.}, hence we deduce that $L\ge 2M$ and that $M$ does not contain any proper subline  bundle $M_1$ with two linearly independent global sections (otherwise\footnote{Write $N=N(S,L)=M_1+K$ with $K$ effective, we get $N\cdot (L-M_1-K)-M_1\cdot (L-M_1)=K\cdot (L-M_1)-N\cdot K=K\cdot (L-M_1-N)>K\cdot (L-2N)\ge 0.$ Recall that $L\ge 2N$.} $N(S,L)=M_1$ ). Hence either $M$ cannot be expressed as the sum of two positive line bundles or $M=K_1+K_2$ with $K_i^2=2$.  Moreover, observe that $c_1(E)=L\ge 2M$, so whenever we have a nonzero map $E\to M$ the kernel is a subline bundle $\ge M$ carrying two linearly independent global sections.

Let us distinguish two cases.
\begin{itemize}
\item $L=2M$ and $\mathrm{NS}(S)\simeq \mathbb Z $. In this case we get by Bogomolov inequality that $F$ is unstable. Indeed $c_2(F)=c_2(E)\le g-g(S,L)+2=2M^2+1-M^2=M^2+1<\frac{1}{3}c_1^2(F)=\frac{4}{3}M^2$ (true as soon as $M^2\ge 4$, this is the case since $M=N(S,L)$ has at least $2$ linearly independent global sections). Hence $F$ is either destabilized by a rank $2$ vector bundle $G$ with $c_1(G)\ge 2M$ or it is destabilized by a subline bundle $\ge M$. In the second situation the induced map $M\to F\to E$ is nonzero (otherwise $M\subset \mathrm{ker}(F\to E)\simeq \mathcal O_S$, contradicting $h^0(M)\ge 2$) hence we get the desired subline bundle  $M\subset E$ with $2$ linearly independent global sections. In the first case up to replacing $G$ with a more positive subsheaf we may suppose the quotient $F/G$ torsion free, we get $c_1(F/G)\le 0$. Now if the induced map $\mathcal O_S\to F/G$ is zero we get a morphism $E\to F/G$ whose kernel is a subline bundle with $c_1\ge 2M$ as desired. If the induced map $\mathcal O_S\to F/G$ is nonzero we get $F/G=\mathcal O_S$ and the sequence defining $F$ is split, contradiction.

\item $L\ge 3M$. Recall that in our range $c_2(E)=d\le g-g(S,L)+2= \frac{(L-M)^2}{2} +\frac{M^2}{2}+1$ and $c_1(E)=L$. A general section $s$ of $E$ induces a sequence $0\to \mathcal O_S \to E \to L \otimes \mathcal I_\xi \to 0$, for $\xi=Z(s)$ of length $c_2(E)$. Twisting by $M^\vee$ we get by standard computations  \begin{eqnarray*}
    \chi(E\otimes M^\vee)&=&\chi(M^\vee)+\chi(L\otimes M^\vee \otimes \mathcal I_\xi) \\
    &=& \frac{M^2}{2}+\frac{(L-M)^2}{2}-c_2(E) \\
    &\ge & M^2-L\cdot M +\frac{L^2}{2}-(\frac{L^2}{2}+1-M\cdot (L-M)) \\
    &=&-1
\end{eqnarray*}
hence
  \[
  \chi(F\otimes M^\vee\otimes \eta^\vee)=\chi(E\otimes M^\vee)+\chi(M^\vee)\ge \frac{M^2}{2}-1>0
  \]
for any $\eta \in \mathrm{Pic}^0(S)$. 
As a consequence, we have a nonzero map, either $F \to M$ or $M \to F$. In the latter case, notice that the image of the map $M\to F$ cannot be contained in $\mathrm{ker}(F\to E)$, hence we obtain a nonzero map $M \to E$. Meanwhile, for the former case hold $h^0(F\otimes M^\vee)=0$ and $h^0(F^\vee \otimes M)=h^2(F\otimes M^\vee)>0$. 
We may suppose the composition $\mathcal{O}_S\to F\to M$ to be nonzero, otherwise we would have a map $E \to M$, thus $E$ would have a subline bundle numerically $\ge 2M$ with two linearly independent global sections. If we had $c_1(\mathrm{Im}(F \to M))=0$, we would get a splitting of the sequence defining $F$, therefore $c_1(\mathrm{Im}(F \to M))>0$. We claim to obtain the exact sequence
  \[
  \begin{tikzcd}
    0 \arrow{r} & G \arrow{r} & F \arrow{r} & M \otimes \mathcal{I}_\xi \arrow{r} & 0
    \end{tikzcd}
  \]
for some $\xi$ finite scheme. If this does not happen, $c_1(\mathrm{Im}(F \to M))>0$ implies that $M=K_1+K_2$ and that the morphism factors through $K_i\otimes \eta \subset M$ (the only positive subline bundles of $M$). Here $M^2\ge 6$ and $\chi(F\otimes M^\vee\otimes \eta)>1$, so we get $H^0(F\otimes M^\vee)$ at least $2-$dimensional. Therefore the claim is equivalent to be able to choose a morphism $F\to M$ which does not factor through $K_i\otimes \eta\subset M$ for some $\eta \in \mathrm{Pic}^0(S)$. First, notice that $\eta$ must be constant, otherwise we would get a nonconstant map $\mathbb P H^0(F\otimes M^\vee)\to \mathrm{Pic}^{K_i}(S)$ which is a contradiction, hence the general morphism $F\to M$ has to factor through $K_i\otimes \eta$ for a fixed $\eta$. Moreover, 
the exact sequence (obtained by taking cohomology in an appropriate twist of the sequence defining $F$)\[ \begin{tikzcd}\mathbb C\simeq H^2(K_i^\vee\otimes \eta^\vee)\to H^2(F\otimes K_i^\vee\otimes \eta^\vee)\to H^2(E\otimes K_i^\vee \otimes \eta^\vee)=0 \end{tikzcd} \] 
yields $h^0(F\otimes \eta^\vee \otimes K_i^\vee)=1$.
The last term of the sequence is supposed to be zero, otherwise $\mathrm{ker}(E\to K_i\otimes \eta)$ would be the desired subline bundle with two linearly independent sections. As claimed, the map $H^0(F^\vee\otimes K_i \otimes \eta)\to H^0(F^\vee\otimes M)$ cannot be surjective and we can find a morphism $F\to M$ surjective outside a finite scheme.

Now pick $\eta_1\neq \eta_2\in \mathrm{Pic}^0(S)$, repeating the construction we get

\[
  \begin{tikzcd}
    0 \arrow{r} & G_i \arrow{r} & F \arrow{r} & M \otimes \eta_i \otimes \mathcal{I}_{\xi_i} \arrow{r} & 0.
    \end{tikzcd}
  \]
for $i=1,2$. Now either $G_1\to M\otimes \eta_2$ or $G_2\to M\otimes \eta_1$ must be nonzero; indeed if the first is zero we get $G_1\subset G_2$ and if the second is zero we get $G_2\subset G_1$. Hence if they both are zero $G_1=G_2$ as subsheaves of $F$, contradiction since the cokernels of the inclusions $G_i\subset F$ are non isomorphic ($\eta_1\neq \eta_2$). The kernel of $G_1\to M\otimes \eta_2$ contains $(L-2M)\otimes \eta_2^\vee$. Hence we get a non zero map \[
(L-M)\otimes \eta_2^\vee\to G_1\to F\to E,
\]
(notice of course that $(L-2M)\otimes \eta_2^\vee$ cannot be contained in $\mathrm{ker}(F\to E)$). Since $L\ge 3M$ we get that $L-2M\ge M$ has at least two linearly independent global sections. 
  
\end{itemize}

\subsubsection{A space parametrizing Lazarsfeld--Mukai bundles}
Now we know we may suppose $E \in
\mathrm{Ext}^1(N \otimes
\mathcal{I}_\xi,M)$ with $h^0(M)\ge 2$. Furthermore, we have the following:
\begin{lemma}
  Let $E,M,N,\xi$ be as above. Then,
  \begin{itemize}
  \item $h^0(N \otimes \mathcal{I}_\xi) \ge 2$ and $|N \otimes \mathcal{I}_\xi|$ has a zero-dimensional base locus.
  
    \end{itemize}
  \label{xl}
\end{lemma}
\begin{proof}
  $E$ is globally generated away from a finite set. This follows from $h^0(E) \ge 3$ together with the exact sequence
  \[
  \begin{tikzcd}
    0 \arrow{r} & H^0(A)^\vee \arrow{r} & H^0(E) \arrow{r} & H^0(\omega_C \otimes A^\vee).
    \end{tikzcd}
  \]
  Indeed since $h^0(E)\ge 3$, the map $H^0(E)\to H^0(\omega_C\otimes A^\vee)$ is non zero. Consider $s\in H^0(E)$ which is non zero when mapped to $H^0(\omega_C\otimes A^\vee)$. We get a commutative diagram with exact rows\[
  \begin{tikzcd}
  (H^0(A)^\vee\oplus \langle s \rangle)\otimes \mathcal O_S \arrow{r} \arrow{d} & E \arrow{r} \arrow{d} & Q\arrow{d}{\simeq} \arrow{r} & 0 \\
  \langle s \rangle \otimes \mathcal O_S \arrow{r} & \omega_C\otimes A^\vee \arrow{r} & \mathcal O_{Z_C(s)} \arrow{r} & 0,
  \end{tikzcd}
  \]
  where in the last term $s$ is seen as a section of $\omega_C\otimes A^\vee$ (hence $Z_C(s)\subset C$). Now since $C$ is a smooth irreducible curve and $\omega_C\otimes A^\vee$ is a line bundle over $C$ we deduce that $Z_C(s)$ is a finite (curvilinear) scheme. This implies that $E$ is globally generated outside $Z_C(s)$.
This implies that also $N \otimes \mathcal{I}_\xi$ must be globally generated away from a finite set (it is a quotient of E).
  
  We immediately get $h^0(N \otimes \mathcal{I}_\xi) \ge 2$ since  $|H^0(N \otimes \mathcal{I}_\xi)|$ has zero dimensional base locus.
\end{proof}
Now we want to reconstruct the general $(C,A)$ starting from an extension with the right properties. For any $[n] \in \mathrm{NS}(S)$ such that $[n]^2>2,([L]-[n])^2>2$,
 we consider
\begin{eqnarray}
{\mathcal{P}}_{[n],l}=\{ (E,N,\eta,\xi)| \quad
N &\in& \mathrm{Pic}^{[n]}(S), \quad \eta \in \mathrm{Pic}^0(S), 
\quad  \xi \in \mathrm{Hilb}^l(S), \nonumber \\   E &\in& \mathrm{Ext}^1(N\otimes \mathcal{I}_{\xi},L \otimes N^\vee \otimes \eta) \textrm{ is a vector bundle}\}. \nonumber 
\end{eqnarray}

Moreover, we consider the relative Grassmanian ${\mathcal{G}}_{[n],l} \to {\mathcal{P}}_{[n],l}$, where  \[({\mathcal{G}}_{[n],l})_{(N,\eta,\xi,E)}=\mathrm{Gr}(2,H^0(E))
.\] 
From a couple $(\Lambda,(N,\eta,\xi,E))$ we may construct a couple $(C,A) \in \mathcal{W}^1_d|L|_{num}$ via the cokernel of the map $E^\vee \to \Lambda^\vee \otimes \mathcal{O}_S$, i.e. we have an exact sequence
\[
\begin{tikzcd}
0 \arrow{r} & E^\vee \arrow{r} & \Lambda^\vee \otimes \mathcal{O}_S \arrow{r} & A \arrow{r} & 0.
  \end{tikzcd}
\]
This gives rise to a rational map $\mathcal{G}_{[n],l} \dashrightarrow \mathcal{W}^1_d|L|_{num}$, this map is not defined in the locus where the cokernel $A$ is not locally free (supported on a singular curve) and in the locus where $\Lambda \otimes \mathcal O_S \to E$ has generic rank $1$ (i.e. it factors through a subline bundle). In principle there may be components of the above Grassmanian where the map is not well defined, we will just remove them. The contents of section 2.1 may be rephrased by saying that any dominating component $\mathcal{W} \subset \mathcal{W}^1_d|L|_{num}$ of dimension $>g+\rho(g,1,d)$ is contained generically in the image of 
\[
 \bigcup_{[n] \in \mathrm{NS}(S), l \in \mathbb{N}} \mathcal{G}_{[n],l} \dashrightarrow \mathcal{W}^1_d|L|.
\]
We will give the required bound on the dimension of $\mathcal{W}$ by studying the dimension of $\mathcal{G}_{[n],l}$ and of the fibers of the above map.

Since $\mathcal{G}_{[n],l} \to \mathcal{P}_{[n],l}$ is non-flat (in particular, the dimension of the fibers becomes bigger as $h^0(E)$ grows), we need a good understanding of the locus where the dimension of the fibers jumps (in particular, we want to understand the dimension).

 We will always compute the dimension locally at smooth points of the respective space. When we will write general we will always mean a general smooth point in its irreducible component. In the following lemma we list subschemes (possibly empty) covering $\mathcal{P}_{[n],l}$, based on the numerical characters of $E$; we postpone a proof.

\begin{lemma}
  Write
    \[
  \mathcal{P}_{[n],l}= \bigcup_{i=0}^4 \mathcal{Z}^i_{[n],l},
  \]
  where
  \begin{eqnarray}
    \mathcal{Z}^i_{[n],l}=\{ E \in \mathcal{P}_{[n],l} \textrm{ so that } h^1(E)=i \}. \nonumber
    \end{eqnarray}
    Then for general $E \in  \mathcal{Z}^i_{[n],l}$ the image of the projection $q:(\mathcal{Z}^i_{[n],l})_{(N,\eta)} \to \mathrm{Hilb}^l(S)$ is at most $(2l-2i)$-dimensional at $q(E)$.
  \label{diff}
\end{lemma}
\begin{proof}
  See section 4.
  \end{proof}
  \subsection{The dimension of $\mathcal{G}_{[n],l}$}
  In the following we write $\mathcal{G}_{[n],l}=\mathcal{G}, \mathcal{P}_{[n],l}=\mathcal{P}$.
  We follow a similar pattern of \cite{FA}, we invite the reader to compare with section 3 of the aforementioned paper.
  In the following, $r$ will denote the rank of the coboundary map $\delta: H^1(\mathcal{O}_S) \to \mathrm{Ext}^2(N \otimes \mathcal{I}_\xi,M)$ (see also the first lines of the following lemma) for $E$ in $\mathcal{P}_{[n],l}$ general in its irreducible component. Let us remark that if the general $E$ has a non trivial endomorphism we can always realize it as $E\in \mathrm{Ext}^1(\mathrm{Im}(\varphi),\mathrm{ker}(\varphi))=\mathrm{Ext}^1(N\otimes \mathcal I_\xi,M)$ for a suitable \footnote{If there is a nontrivial endomorphism $\varphi$ we can find a suitable linear combination between $\varphi$ and $\mathrm{id}$ dropping rank at a point. Since the determinant of an endomorphism is a section of $\mathcal O_S$ we get that the chosen linear combination drops rank everywhere. } $\varphi$ (see also \cite{FA}); we will always make such a choice if possible. \\
  \begin{lemma}
   Consider a general $E \in \mathrm{Ext}^1(N \otimes \mathcal{I}_{\xi},M)$ and recall $l=\mathrm{length}(\xi)$, we have:
    \begin{enumerate}
      \item 
     $ \mathrm{dim}(\mathrm{Ext}^1(N \otimes \mathcal{I}_{\xi},M)) \le  l+h^1(M \otimes N^\vee);$

  \item 
      $\mathrm{dim}(\mathrm{Ext}^1(N \otimes \mathcal{I}_{\xi},M)) \le  r-1+h^1(E^\vee \otimes M).$
    
    \end{enumerate}
Suppose now that $E\in \mathcal{P}$ is a Lazarsfeld--Mukai bundle  with non-trivial endomorphisms $h^0(E \otimes E^\vee) >1$ write $E\in \mathrm{Ext}^1(\mathrm{Im}(\varphi),\mathrm{ker}(\varphi))=\mathrm{Ext}^1(N\otimes \mathcal I_\xi,M)$. We have:
    \begin{eqnarray}
      h^0(E \otimes E^\vee)&=&1+h^0(M \otimes N^\vee) \quad \textrm{if} \quad l>0 \nonumber, \\
      h^0(E \otimes E^\vee)&=&2+h^0(M \otimes N^\vee)+h^0(M^\vee \otimes N) \quad \textrm{otherwise.} \nonumber
    \end{eqnarray}
   
    \label{label}
    \end{lemma}

  \begin{proof}
    For the statements about the endomorphisms see \cite{FA}.
    For the first part, applying the functor $\mathrm{Hom}_S(-,M)$ to the  exact sequence
    \[
    \begin{tikzcd}
      0 \arrow{r} & M \arrow{r} & E \arrow{r} & N \otimes \mathcal{I}_\xi \arrow{r} & 0,
    \end{tikzcd}
    \]
    we get in cohomology
    \[
    \begin{tikzcd}
      \mathbb{C} \arrow{r} & \mathrm{Ext}^1(N\otimes \mathcal{I}_\xi, M) \arrow{r} & H^1(E^\vee \otimes M) \arrow{r} & H^1(\mathcal{O}_S) \arrow{r}{\delta} & \mathrm{Ext}^2(N \otimes \mathcal{I}_\xi, M).
      \end{tikzcd}
    \]
    We deduce $(2)$.
    
   For $(1)$ notice that $\mathrm{Ext}^1(N\otimes \mathcal I_\xi, M)=\mathrm{Ext}^1(N\otimes M^\vee \otimes \mathcal I_\xi, \mathcal O_S)=H^1(N\otimes \mathcal I_\xi\otimes M^\vee)$, the last equality by Serre duality. Moreover we have an exact sequence $0\to N\otimes M^\vee \otimes \mathcal I_\xi \to N \otimes M^\vee \to \mathcal O_\xi \to 0,$ inducing in cohomology $H^0(\mathcal O_\xi)\to H^1(N\otimes M^\vee \otimes \mathcal I_\xi)\to H^1(M\otimes N^\vee)$; $(1)$ follows.
  \end{proof}
   Let $\pi:\mathcal{G} \to \mathcal{P}$ be the projection. In the notation of lemma \ref{diff}, we can write
    \[
    \mathcal{G}=\bigcup_{i=0}^4 {\overline{\pi^{-1}(\mathcal{Z}^i)}}.
    \]
      We can prove the following:
    \begin{theorem}
      If $h^2(M \otimes N^\vee)=0$ or if $M\otimes N^\vee=\mathcal{O}_S$ for the general $E$ then
      \begin{eqnarray}
      \mathrm{dim}(\mathcal{G}) &\le& g+l+h^0(M \otimes N^\vee)-2 \quad \quad \textrm{if } \quad l>0; \nonumber \\
      \mathrm{dim}(\mathcal{G}) &\le& g+h^0(M \otimes N^\vee)-1 \quad \quad \textrm{if } \quad l=0. \nonumber
      \end{eqnarray}
\label{stoca}
    \end{theorem}
    \begin{proof}
    We divide the proof in several cases.
    \begin{itemize}
        \item
    Suppose $l=0$ and $h^2(M \otimes N^\vee)=0$. We get $E\in \mathrm{Ext}^1(N,M)$ and $h^1(E)=0$. Hence $\mathrm{dim}(\mathcal{P}) \le \mathrm{dim}(\mathrm{Pic}^{[n]}(S))+\mathrm{dim}(\mathrm{Pic}^0(S))+h^1(M\otimes N^\vee) =4+h^1(M\otimes N^\vee)$, with equality holding if and only if $h^1(M\otimes N^\vee)=0$. Now $\mathcal{G}$ is a $\mathrm{Gr}(2,H^0(E))-$bundle over $\mathcal{P}$, hence it is of relative dimension $2(h^0(E)-2)=2(g-d-3)$. By Riemann-Roch $h^0(M\otimes N^\vee)=\chi(M \otimes N^\vee)+h^1(M\otimes N^\vee)=\chi(M\otimes N)-2M\cdot N+h^1(M\otimes N^\vee)=g-1-2M\cdot N+h^1(M\otimes N^\vee)=g-1-2d+h^1(M\otimes N^\vee)$. Putting all together we get \[
    \mathrm{dim}(\mathcal{G})=\mathrm{dim}(\mathcal P)+2(h^0(E)-2)\le 4+h^1(M\otimes N^\vee)+2g-2d-6=g-1+h^0(M\otimes N^\vee),
    \]
with equality holding if $h^1(M\otimes N^\vee)=0$.
  \item  Suppose $l=0$ and $M=N$. We get $E\in \mathrm{Ext}^1(M,M)$ and $h^1(E)=0$. Hence $\mathrm{dim}(\mathcal{P}) = \mathrm{dim}(\mathrm{Pic}^{[n]}(S))+h^1(M\otimes N^\vee)-1 =3$ (the bundle $E$ is determined by $M$ and the extension class in $\mathbb{P}(\mathrm{Ext}^1(M,M))$ which is $1-$dimensional). Now $\mathcal{G}$ is a $\mathrm{Gr}(2,H^0(E))-$bundle over $\mathcal{P}$, hence it is of relative dimension $2(h^0(E)-2)=2(g-d-3)$. By Riemann-Roch $h^0(M\otimes N^\vee)=\chi(M \otimes N^\vee)+h^1(M\otimes N^\vee)-h^2(M\otimes N^\vee)=\chi(M\otimes N)-2M\cdot N+h^1(M\otimes N^\vee)-h^2(M\otimes N^\vee)=g-1-2M\cdot N+h^1(M\otimes N^\vee)-1=g-2-2d+h^1(M\otimes N^\vee)$. Putting all together we get \[
    \mathrm{dim}(\mathcal{G})=\mathrm{dim}(\mathcal P)+2(h^0(E)-2)= 1+h^1(M\otimes N^\vee)+2g-2d-6<g-1+h^0(M\otimes N^\vee),
    \]
hence the Theorem.

    \item  Suppose $l>0$ and $h^2(M\otimes N^\vee)=0$.  It suffices to prove the inequality for the components $\mathcal{G}=\bigcup_{i=0}^4 {\overline{\pi^{-1}(\mathcal{Z}^i)}}$. 
        We can write
        \[
        \mathrm{dim}(\pi^{-1}(\mathcal{Z}^i))\le 2(l-i)+\mathrm{dim} \mathbb{P}(\mathrm{Ext}^1(N \otimes \mathcal{I}_{\xi},M))+\mathrm{dim}[n]+\mathrm{dim}(\mathrm{Pic}^0(S))+2(h^0(E)-2);
        \]
        where the first term is an upper bound of the dimension of the scheme where $\xi$ moves and the last one is an upper bound of the dimension of the general fiber \footnote{Recall that $\mathcal{G}\to \mathcal{P}$ is a a relative $\mathrm{Gr}(2,H^0(E))$ bundle. Hence the fibers have dimension $\mathrm{dim}(\mathrm{Gr}(2,H^0(E))=2(h^0(E)-2)$.} of $\pi_{|_{\pi^{-1}(\mathcal{A}^i)}}:\mathcal{G} \to \mathcal{P}$. Now, recall that $h^0(E)=g-1-d+i$ and $M\cdot N+l=c_2(E)=d$. By Riemann-Roch $h^0(M\otimes N^\vee)-h^1(M\otimes N^\vee)=\chi(M \otimes N^\vee)=\chi(M\otimes N)-2M\cdot N=g-1-2M\cdot N$.
        Putting all together we get:
        \begin{eqnarray}
          \mathrm{dim}(\pi^{-1}(\mathcal{Z}^i))&\le& 3d -3M\cdot N+h^1(M \otimes N^\vee)+2g-3-2d \nonumber \\
          &=& g+d-M\cdot N+h^0(M\otimes N^\vee)-2 \nonumber \\
          &=& g+l+h^0(M \otimes N^\vee)-2. \nonumber
        \end{eqnarray}
\item Suppose $l>0$ and $M=N$.  It suffices to prove the inequality for the components $\mathcal{G}=\bigcup_{i=0}^4 {\overline{\pi^{-1}(\mathcal{Z}^i)}}$. 
        We can write
        \[
        \mathrm{dim}(\pi^{-1}(\mathcal{Z}^i))\le 2(l-i)+\mathrm{dim} \mathbb{P}(\mathrm{Ext}^1(M \otimes \mathcal{I}_{\xi},M))+\mathrm{dim}[n]+2(h^0(E)-2);
        \]
        where the first term is an upper bound of the dimension of the scheme where $\xi$ moves and the last one is an upper bound of the dimension of the general fiber of $\pi_{|_{\pi^{-1}(\mathcal{A}^i)}}:\mathcal{G} \to \mathcal{P}$. Now, recall that $h^0(E)=g-1-d+i$ and $M\cdot N+l=c_2(E)=d$. By Riemann-Roch $h^0(M\otimes N^\vee)-h^1(M\otimes N^\vee)=\chi(M \otimes N^\vee)-1=\chi(M\otimes N)-2M\cdot N-1=g-2-2M\cdot N$.
        Putting all together we get:
        \begin{eqnarray}
          \mathrm{dim}(\pi^{-1}(\mathcal{Z}^i))&\le& 3d -3M\cdot N+h^1(M \otimes N^\vee)+2g-3-2d \nonumber \\
          &<& g+d-M\cdot N+h^0(M\otimes N^\vee)-2 \nonumber \\
          &=& g+l+h^0(M \otimes N^\vee)-2. \nonumber
        \end{eqnarray}

         \end{itemize}
      \end{proof}
In the following section we will need the following technical lemma on the projection $\mathcal{P}_{[n],l} \to \mathrm{Pic}^0(S)$ in terms of $\Psi_{E,M,N}=\{\eta \in \mathrm{Pic}^0(S) | \quad h^2(M \otimes \eta \otimes E^\vee)>0\}$. 
    \begin{lemma}
    Suppose $r>0$ and $M \neq N$ for $E$ general in its irreducible component. Moreover, suppose $\Psi_{E,M,N}$ is of dimension $i$ and  fix $N,\xi$. Then, for a sufficiently small open neighborhood of $E \in \mathcal{P}$, the image of the projection $U\cap  (\mathcal{P}_{[n],l})_{\xi,N} \to  \mathrm{Pic}^0(S)$ is  at most $i-$dimensional.
        \label{lastp}    \end{lemma}
    \begin{proof}
  Suppose $r >0$ for some $E$. Take $\eta$ in a neighborhood of $\mathcal{O}_S$ and apply the functor $\mathrm{Hom}(-,M\otimes \eta)$ to the same sequence as before. We get
     \[
    \begin{tikzcd}
      0=H^0(\eta) \arrow{r} & \mathrm{Ext}^1(N\otimes \mathcal{I}_\xi, M\otimes \eta) \arrow{r} & H^1(E^\vee \otimes M \otimes \eta) \arrow{r} & H^1(\eta)=0.
      \end{tikzcd}
    \]
Now saying that that $\eta\in \mathrm{Pic}^0(S)$ is in the image of the projection of a small analytic neighborhood $U$ of $E$ in $\mathcal P$ is equivalent to say that $E$ can be deformed as a Lazarsfeld--Mukai bundle in $\mathrm{Ext}^1(N \otimes \mathcal{I}_\xi, M \otimes \eta)$. Now, supposing $E$ and $\eta $ general (we are always supposing $E$ general in its irreducible component\footnote{In other words we are choosing $E$ general in its irreducible component such that the value of $\mathrm{dim}(\mathrm{Ext}^1(N \otimes \mathcal{I}_\xi, M))$ is minimal among the Lazarsfeld--Mukai bundles in the same irreducible component. $\eta$ is general in the image of the projection. }) we must have $\mathrm{dim}(\mathrm{Ext}^1(N \otimes \mathcal{I}_\xi, M \otimes \eta))=\mathrm{dim}(\mathrm{Ext}^1(N \otimes \mathcal{I}_\xi, M))$. We get
    \[
    h^1(E^\vee \otimes M \otimes \eta)=h^1(E^\vee \otimes M)-r+1.
    \]
    Now, since $h^0(E^\vee \otimes M \otimes \eta)=h^0(E^\vee \otimes M)=0$ and $\chi(E^\vee \otimes M)=\chi(E^\vee \otimes M \otimes \eta)$, we get $h^2(E^\vee \otimes M \otimes \eta)=h^2(E^\vee \otimes M)+r-1>0$. The facts concerning the projection $U\cap  (\mathcal{P}_{[n],l})_{\xi,N} \to  \mathrm{Pic}^0(S)$ follow: we just proved that if we can deform $E$ as a Lazarsfeld--Mukai bundle in $\mathrm{Ext}^1(N \otimes \mathcal{I}_\xi, M \otimes \eta)$ (hence if $\eta $ is in the image of the projection)  then $h^2(E^\vee \otimes M\otimes \eta)>0$.
   
    \end{proof}

    \subsection{A bound for the dimension of $\mathcal{W}^1_d|L|$}
   Let us recall the rational map $h_{[n],l}:\mathcal{G}_{[n],l} \dashrightarrow \mathcal{W}^1_d(|L|_{num})$, where $h_{[N,l]}(E,\Lambda)=(C,A)$; the couple $(C,A)$  is induced via the following exact sequence:
    \[
    \begin{tikzcd}
      0 \arrow{r} & \Lambda \otimes \mathcal{O}_S \arrow{r} & E \arrow{r} & \omega_C \otimes A^\vee \arrow{r} & 0.
      \end{tikzcd}
    \]
   
    We have the following:
    \begin{lemma}
    If $l>0$ the fibers of the map $h_{[n],l}$ have dimension $\ge h^0(M\otimes N^\vee)$, if $l=0$ the fibers of the map $h_{[n],l}$ have dimension $\ge h^0(M\otimes N^\vee)+1$.
    \label{izzato}
    \end{lemma}
    \begin{proof}

     The fiber of the $h_{[n],l}$ contains a Zariski open subset of $\mathbb{P}H^0(\mathrm{End}(E))$. Indeed let $\varphi :E \to \omega_C \otimes A^\vee$ be the canonical morphism, then for any automorphism $\iota \neq id$ of $E$ we get $\varphi\iota \neq \varphi$ as the following argument shows.
      Applying $-\otimes E^\vee$ to the sequence defining $E=E_{C,A}$ we get
      \[
      \begin{tikzcd}
        0 \arrow{r} & H^0(A)^\vee \otimes E^\vee \arrow{r} & E^\vee \otimes E \arrow{r} & E^{\vee} \otimes \omega_C \otimes A^\vee \arrow{r} & 0.
        \end{tikzcd}
      \]
      Hence in cohomology
      \[
      \begin{tikzcd}
        0 \arrow{r} & H^0(\mathrm{End}(E)) \arrow{r} &H^0(\mathrm{Hom}(E, \omega_C \otimes A^\vee)) \arrow{r} & H^0(A)^\vee \otimes H^1(E).
        \end{tikzcd}
      \]
      By lemma \ref{label}, if $l>0$
      \[
      h^0(\mathrm{Hom}(E, \omega_C \otimes A^\vee)) \ge h^0(\mathrm{End}(E))=1+h^0(M\otimes N^\vee).
      \]
      hence the claim.  The claim for $l=0$ is analogous.
    \end{proof}
    Putting it all together we get:
    \begin{lemma}
      $\mathrm{dim}(\mathrm{Im}(h_{[n],l})) \le g+l-2=g+d-M\cdot N-2$.
      \label{oie}
    \end{lemma}
    \begin{proof}
      The case $h^2(M \otimes N^\vee)=0$ or $M=N$ follows from Theorem \ref{stoca} together with lemma \ref{izzato}. We shall still cover the case when $h^2(M \otimes N^\vee)>0$ and $M \neq N$ for the general $E$. If $r=0$ the fibers of the map $h_{[n],l}$ contain $\mathbb{P}H^0(E \otimes M^\vee)=\mathbb{P}H^2(E^\vee \otimes M)^\vee$, if $r>0$ the fibers of the map contain $\bigcup_{\eta \in \mathrm{Pic}^0(S)}\mathbb{P}H^0(E \otimes M^\vee \otimes \eta)$ (see also lemma \ref{lastp}). Indeed for any morphism $f:M\otimes \eta \to E$ we realize $E$ as an element of $\mathrm{Ext}^1(N'\otimes \mathcal{I}_{\xi'},M\otimes \eta)$ where $N' \otimes \mathcal{I}_{\xi'}=\mathrm{coker}(f)$. Moreover, all these extensions are non-isomorphic since we can suppose that $E$ has trivial endomorphism \footnote{If $E$ has an endomorphism and we realize it in $\mathrm{Ext}^1(\mathrm{Im}(\varphi),\mathrm{ker}(\varphi))$ we may always suppose $h^2( M \otimes N^\vee)=0$}. We cover just the case $r>0$ and $\Psi_{E,M,N}=\mathrm{Pic}^0(S)$ (i.e. $h^2(E^\vee \otimes M \otimes \eta)>0$ for any $\eta$), the other cases being similar.  We get (analogously to Theorem \ref{stoca}, the last term below coming from the dimension of the fibers): 
      \begin{eqnarray}
        \mathrm{dim} (\mathrm{Im} h_{[n],l})& \le & 2(l-i)+ \mathrm{dim}(\mathrm{Ext}^1(N \otimes \mathcal{I}_\xi,M))+4-1+2(h^0(E)-2)-(h^2(E^\vee \otimes M) \nonumber+1) \\
        &=& 2(l-i)+2+2(g-1-d+i-2)-\chi(E^\vee \otimes M)-3 \nonumber
        \\
        &=&g+l-2. \nonumber
      \end{eqnarray}
   
      \end{proof}

 Now, $M\cdot N \ge g(S,L)-2$ by definition, so Theorem \ref{t3} follows easily.

In order to conclude the proof of Theorem \ref{t1}, we need a lower bound on the gonality.
\begin{theorem}
  Let $S$ be a simple abelian surface. Let $L=|C|$ be an effective linear system on $S$, then $\mathrm{gon}(C) \ge g(S,L)-2$ (see also Theorem \ref{t1}).
\end{theorem}
\begin{proof}
  Let $A \in W^1_d(C)$ be a line bundle computing the gonality, let us consider the associated Lazarsfeld--Mukai bundle $E=E_{C,A}$. If $E$ is simple and stable, we have by Bogomolov inequality
  \[
  \mathrm{gon}(C)=d=c_2(E) \ge \frac{1}{4}c_1(E)^2=\frac{1}{2}(g-1)= \frac{1}{2}(g+3)-2 \ge g(S,L)-2.
  \]
  If $E$ is not simple and stable, it comes equipped with a short exact sequence
  \[
  \begin{tikzcd}
    0 \arrow{r} & M \arrow{r} & E \arrow{r} & N \otimes \mathcal{I}_\xi \arrow{r} &0.
    \end{tikzcd}
  \]
  In this case, we get $d=c_2(E)=M\cdot N+\mathrm{length}(\xi) \ge g(S,L)-2$. Indeed $h^0(N)\ge 2$ because it is globally generated outside a finite set and $h^0(M)\ge 2 $ because $M \to E$ destabilizes $E$, so we may suppose $M^2 \ge N^2$.
  \end{proof}
\section{An upper bound for the gonality and linear growth condition}
In this section, we want to give an upper bound for the gonality of any curve on an abelian surface. It will turn out as in the case of $K3$ surfaces that for general curves in the linear systems the special behaviour of the curves is inherited from the surface, but unlike the $K3$ case, the Clifford index is not computed by a line bundle coming from the surface. The idea is that the map $C \to \mathbb{P}^1$ computing the gonality is the restriction of a (rational) map $\varphi: S \dashrightarrow \mathbb{P}^1$ obtained by imposing the maximal possible number of base points on $C$. In the following $N=N(S,L)$ (see the introduction for a definition).  Let us first deal with the case $h^0(N)=2$ which is more elementary. Notice that in case $L=2N$ and $\mathrm{NS}(S)=\mathbb Z \cdot N$ we do not need to prove anything since we showed that the curves are of maximal gonality.
\begin{lemma}
  Given any two points $P,Q \in S$ and any curve $C \in |L|$ there exists $a \in S$ such that $P,Q \in  C+a$.
\end{lemma}
\begin{proof}
  Up to translation, we may suppose that $P=0$ the origin in $S$. Let us consider the map
  \[
  \begin{tikzcd}
    C \times C \arrow{r}{f} & S,
\end{tikzcd}
  \]
  where     $f(c_1,c_2)   = c_1-c_2$. Since $C$ is ample the map $f$ is surjective, hence $Q=d_1-d_2$ for some $d_1,d_2 \in C$. Take $a=-d_2$, you get $P=0=d_2-d_2 \in C+a$ and $Q=d_1-d_2 \in C+a$. 
\end{proof}
\begin{corollary}
  Suppose $h^0(N)=2$, then each curve in $L$ carries at least six $g^1_{g(S,L)}$. In particular, Theorem \ref{t1} holds for $|L|$. 
\end{corollary}
\begin{proof}
  The rational map $\varphi_{|_N}:S \dashrightarrow \mathbb{P}^1$ has four base points. Apply the previous lemma to any couple of them and you get a pencil of degree $g(S,L)$ on $C+a$.
  \end{proof}

Now we deal with the general case. Let us settle some notation.

Let us consider the Grassmaniann of planes $\mathcal{G} \to \mathrm{Pic}^{[n]}(S)$ (the fiber over $N \otimes \eta$ is \\ $\mathrm{Gr}(2, H^0(N \otimes \eta))$). We have a rational map of finite degree onto its image  $\mathcal{G} \dashrightarrow \mathrm{Sym}^{N^2}S$ given by sending $\langle s_1,s_2 \rangle \in \mathrm{Gr} (2, H^0(\mathcal{O}_S(N) \otimes \eta))$ to $Z(s_1)\cap Z(s_2)$. Let us consider the pullback of this locus to $S^{2h^{0}(N)}$ via the cartesian square
\[
\begin{tikzcd}
  \mathcal{G}' \arrow[dashed]{r} \arrow[dashed]{d}{\gamma} & \mathcal{G} \arrow[dashed]{d} \\
  S^{N^2} \arrow{r} & \mathrm{Sym}^{N^2}(S).
  \end{tikzcd}
\]
Finally let us consider  $F_0,F_1,F_2 \subset \mathcal G'\times H^0(L)$, generically defined as
\begin{eqnarray}
  F_0&=&\{(p,s) \in \mathcal G'\times H^0(L) \quad | \, s\in H^0(L\otimes \mathcal I_{P_1,\dots,P_{N^2}}), \gamma(p)=(P_1,\dots,P_{N^2})\} \nonumber \\
   F_1&=&\{(p,s) \in \mathcal G'\times H^0(L) \quad | \, s\in H^0(L\otimes \mathcal I_{P_1,\dots,P_{N^2-1}}), \gamma(p)=(P_1,\dots,P_{N^2})\} \nonumber \\
   F_2&=&\{(p,s) \in \mathcal G'\times H^0(L) \quad | \, s\in H^0(L\otimes \mathcal I_{P_1,\dots,P_{N^2-2})}, \gamma(p)=(P_1,\dots,P_{N^2})\}. \nonumber
\end{eqnarray}

The construction makes sense birationally (in the locus where the points are distinct; notice that, since we may suppose $N^2>4$, the linear system $|N|$ has no base points).  
We have canonical maps $F_0,F_1,F_2 \to H^0(L)$. The goal of this section will be to prove that these maps are generically finite onto their images, then a quick parameters count yields that $F_2 \to H^0(L)$ is dominant. Equivalently, the general curve in $|L|$ carries at least one $g^1_{g(S,L)}$ (recall $g(S,L)=N\cdot (L-N)+2$).
We first prove that $F_0 \to H^0(L)$ is of finite degree onto its image. We consider $C$ a general smooth curve in the image of $F_0$, this means $C \in H^0(L \otimes \mathcal{I}_{P_1,\dots,P_{N^2}})$ for $P_1,\dots,P_{N^2} \in Z(s_1) \cap Z(s_2)$ ($s_1,s_2 \in H^0(N\otimes \eta)$ for some $\eta \in \mathrm{Pic}^0(S)$). The two sections $s_1,s_2$ generate a pencil $A$ on $C$ of degree $N\cdot (L-N)=g(S,L)-2$. Hence we have an injective morphism $\mathbb{P}(F_0) \to \mathcal{W}^1_{g(S,L)-2}|L|$. We need the following preliminary lemma:
\begin{lemma}
Suppose $L\ge 3N$, then for general $(P_1,\dots,P_{N^2}) \in \mathcal{G}$ and any tangent directions $v_1,\dots ,v_{N^2} \in T_0(S)$ there exists a smooth curve $C \in H^0(L \otimes \mathcal{I}_{P_1,\dots,P_{N^2}})$ with tangent direction $v_i$ at $P_i$.
\label{dir}
\end{lemma}
\begin{proof}
Let us fix general points $P_1,\dots, P_{N^2}\in \mathcal{G}$ and arbitrary tangent directions $v_i \in T_{P_i}S$.
 Since $L-N \cdot N \ge 2 N^2$, by Riemann-Roch applied to $H_1=Z(s_1)$ we get that $P_1, \dots, P_{N^2}$ impose independent conditions on $H^0(L \otimes N^\vee\otimes \eta)$ for   $\eta \neq \mathcal{O}_S \in \mathrm{Pic}^0(S)$. This follows from the exact sequence
 \[
\begin{tikzcd}
0 \arrow{r} & L  \otimes N^{\otimes -2}\otimes \eta \arrow{r} & L \otimes N^{\vee} \otimes \eta \arrow{r} & L \otimes N^\vee\otimes \eta_{|_{H_1}} \arrow{r} & 0.
\end{tikzcd}  
\]
Hence for any $j=1,\dots, N^2$ we may find a section $Q_j \in H^0(L \otimes N ^{\vee} \otimes \eta)$ passing through all the marked points but $P_j$. Moreover, for any $j$ we may choose $l_j \in \langle s_1,s_2 \rangle$ with tangent direction $v_j$ at $P_j$ (we may suppose $s_1,s_2$ smooth with transverse intersection). 
 Now a general linear combination \[
\sum_j \lambda_j l_j \otimes Q_j
\]
will define an element of $H^0(L \otimes \eta \otimes  \mathcal{I}_{P_1,\dots,P_{N^2}})$ smooth at the marked points with the desired tangent directions. This finishes the proof of the lemma.

\end{proof}

\begin{lemma}
Suppose $L\ge 3N$ and $h^0(N)\ge 3$.
  Let $(C,A) \in \mathcal{W}^1_{g(S,L)-2}|L|$ be a general element in the image of $\mathbb{P}(F_0)$, then $h^0(A^{\otimes 2})=3$.
\end{lemma}
\begin{proof}
  Let us remark that $H^0(N)=H^0(N_{|_{C}})$ and $H^0(N^{\otimes 2})=H^0(N^{ \otimes 2}_{|_{C}})$, this follows from $h^1(L \otimes N^{\otimes -2})=0$. So $H^0(A^{\otimes 2})$ are just sections of $H^0(N^{\otimes 2})$ tangents to $C$ at $P_1,\dots,P_{N^2}$, let us call this vector space $V$. The goal is to show that $V=\langle s_1 \otimes s_1,s_1\otimes s_2,s_2 \otimes s_2\rangle$ for general $(C,A)$. We argue by contradiction, suppose that $\mathrm{dim}(V) \ge 4$. First let us observe that the general element of $V$ must be smooth at least at a point, otherwise the map $\varphi_V:S \to \mathbb{P}^3$ is non-degenerate with a curve as image (it is of degree zero on the elements of $|V|$). The degree of the curve must be at least $3$, hence $H_1+H_2$ cannot be a hyperplane section of such a map. Since being smooth at a point is an open condition, for general $C$ the general element of $V$ is smooth at $P_1$ (up to reindexing). Now consider an irreducible component $\mathcal{G}''\subset \mathcal{G}'$ and let $G''$ be the subgroup of the symmetric group fixing such a component, it is easy to observe that $G''$ is transitive if $h^0(N) \ge 3$. Hence the general element of $V$ for general $C$ is smooth at all the marked points. Now by the previous lemma for $P_1,\dots,P_{N^2}$ general we may find a smooth $C$ with any tangent directions, hence if we choose general tangent directions we may find a smooth quadric $Q \in H^0(N^{\otimes 2})$ with generic tangent directions (the non decomposable element of $V$, whose existence is supposed by contradiction). This implies that imposing tangency with $C$ at $P_1,\dots,P_{N^2}$ impose $N^2$ independent conditions on $H^0(N^{\otimes 2} \otimes \mathcal{I}_{P_1,\dots,P_{N^2}})$. The space $H^0(N^{\otimes 2} \otimes \mathcal{I}_{P_1,\dots,P_{N^2}})$ is of codimension $N^2-1$ as a vector subspace of $H^0(N^{\otimes 2})$. Hence we get $\mathrm{dim}(V)=1$ for general $(C,A)$ contradicting our very first hypothesis.
\end{proof}
The above lemma is the key part to describing the tangent space  $T_{A\otimes \mathcal{O}_C(P_{N^2-1}+P_{N^2})}W^1_{g(S,L)}(C)$.
\begin{lemma} 
  For general $(C,A \otimes \mathcal{O}_C(P_{N^2-1}+P_{N^2})) \in \mathcal{W}^1_{g(S,L)}|L|$ in the image of the canonical map $F_0 \subset F_2 \to \mathcal{W}^1_{g(S,L)}|L|$ we have $\mathrm{dim} (T_{A\otimes \mathcal{O}_C(P_{N^2-1}+P_{N^2})}W^1_{g(S,L)}(C))=2$.
\end{lemma}
\begin{proof}

  Recall that $\mathrm{dim}(T_{A\otimes \mathcal{O}_C(P_{N^2-1}+P_{N^2})}W^1_{g(S,L)}(C))=\mathrm{corank}(\mu)$  where $\mu$ is the Petri map
  \[
  \begin{tikzcd}
    \mu: H^0(A\otimes \mathcal{O}_C(P_{N^2-1}+P_{N^2})) \otimes H^0(\omega_C \otimes A^\vee \otimes \mathcal{O}_C(-P_{N^2-1}-P_{N^2})) \arrow{r} & H^0(\omega_C).
  \end{tikzcd}
  \]
  Hence we deduce from the previous lemma that the multiplication map $H^0(A) \otimes H^0(\omega_C \otimes A^\vee) \to H^0(\omega_C)$ is surjective. We want to prove that $\mathrm{corank}(\mu)=2$ for the Petri map associated with $A \otimes \mathcal{O}_C(P_{i}+P_{j})$ for some $i\neq j$. We have a commutative diagram
  \[
  \begin{tikzcd}
    H^0(A(P_{i}+P_j)) \otimes H^0(\omega_C \otimes A^\vee(-P_{i}-P_j)) \arrow{r} \arrow[hook]{d} & H^0(\omega_C(-P_i-P_j)) \arrow[hook]{r} &  H^0(\omega_C) \arrow{d}{=} \\
    H^0(A)\otimes H^0(\omega_C) \arrow{rr} & & H^0(\omega_C).
    \end{tikzcd}
  \]
  Let $K=\mathrm{ker}(H^0(A) \otimes H^0(\omega_C \otimes A^\vee) \to H^0(\omega_C))$, we want to prove that \\$K \cap H^0(A) \otimes H^0(\omega_C \otimes A^\vee \otimes \mathcal{O}_C(-P_i-P_j)) \subset K$ is of codimension $2$ for some $i \neq j$.
  We consider the natural inclusion $H^0(A) \subset H^0(A \otimes \mathcal{O}_C(\sum_{j=1}^{N^2}P_j))$. Now, since $A \otimes \mathcal{O}_C(\sum_{j=1}^{N^2}P_j)=\mathcal{O}_S(N)_{|_C}$, we get a commutative diagram
  \[
  \begin{tikzcd}
    H^0(A) \otimes H^0(L \otimes N^\vee) \arrow{r}\arrow[hook]{d} & H^0(L) \arrow{d} \\
    H^0(A) \otimes H^0(\omega_C \otimes A^\vee) \arrow{r} & H^0(\omega_C),
    \end{tikzcd}
  \]
  inducing an isomorphism between  $\mathrm{ker}(H^0(A) \otimes H^0(L \otimes N^\vee) \to H^0(\omega_C))$ and \[
  K'=K \cap H^0(A) \otimes H^0(\omega_C \otimes A^\vee \otimes \mathcal{O}_C(-\sum_{j=1}^{N^2}P_j))\subset K.
  \] Let us remark that $\mathrm{rk}(H^0(A) \otimes H^0(L \otimes N^\vee) \to H^0(\omega_C))=\mathrm{rk} (H^0(A) \otimes H^0(L \otimes N^\vee) \to H^0(L))-1$ since the section $s \in H^0(L)$ defining $C$ lies in $H^0(L \otimes \mathcal{I}_{P_1,\dots,P_{N^2}})$.
  The kernel of the multiplication $H^0(A) \otimes H^0(L \otimes N^\vee) \to H^0(L)$ is $H^0(L \otimes N^{\otimes -2})$ by the exact sequence
  \[
  \begin{tikzcd}
    0 \arrow{r} & L \otimes N^{\otimes -2} \arrow{r} & H^0(A) \otimes L \otimes N^\vee \arrow{r} & L \otimes \mathcal{I}_{P_1,\dots,P_{N^2}} \arrow{r} & 0.
    \end{tikzcd}
  \]
  We deduce $\mathrm{dim}(K')= h^0(L \otimes N^{\otimes -2})+1=\frac{L^2}{2}-2L\cdot N+2N^2+1=\mathrm{dim}(K)-2$. This implies there exists $i \neq j$ such that $K \cap H^0(A) \otimes H^0(\omega_C \otimes A^\vee \otimes \mathcal{O}_C(-P_i-P_j)) \subset K$
 is of codimension $2$.
  \end{proof}

Now we are finally ready to prove:
\begin{theorem}
  The two morphisms $\mathbb{P}(F_0),\mathbb{P}(F_2) \to |L|$ are generically finite onto their image.
\end{theorem}
\begin{proof}
  The map $\mathbb{P}(F_0) \to |L|$ factors through the morphism $\mathbb{P}(F_0) \to \mathcal{W}^1_{g(S,L)-2}|L|$, which is injective up to the action of the symmetric group permuting the base points (if two maps defined on $S$, $\gamma_1,\gamma_2 :S \to \mathbb{P}^1$ induce the same morphism $C \to \mathbb{P}^1$ then the fibers in $S$ of the two maps above $\infty$ must intersect in $L\cdot  N-N^2>N^2$ points, i.e. the fiber above $\infty$ in $C$).
  Since  $T_AW^1_{g(S,L)-2}(C)=0$, $A$ cannot be deformed in $W^1_{g(S,L)-2}(C)$. Hence the fibers of $W^1_{g(S,L)-2}(C) \to |L|$ are finite schemes locally at points coming from $\mathbb{P}(F_0)$. Hence  $\mathbb{P}(F_0) \to |L|$ is generically finite.

  Now let us prove that  $\pi:\mathbb{P}(F_2) \to |L|$ is generically finite. We have $\mathbb{P}(F_0) \subset \mathbb{P}(F_2)$, we will prove that the fibers of $\pi$ near a general point of $\mathbb{P}(F_0)$ are zero-dimensional.
  The map $\mathbb{P}(F_2) \to |L|$ factors through the generically finite morphism $f_2:\mathbb{P}(F_2) \to \mathcal{W}^1_{g(S,L)}$. It suffices to show that $\mathrm{Im}(f_2)\cap W^1_{g(S,L)}(C)$ is zero dimensional at a point. Let $(C,A\otimes \mathcal{O}_C(P_{N^2-1}+P_{N^2}))$ be a general point in the image of $\mathbb{P}(F_0)$. By the previous lemma $\mathrm{dim}(T_{A\otimes \mathcal{O}_C(P_{N^2-1}+P_{N^2})}W^1_{g(S,L)}(C))=2$. The line bundle $A\otimes \mathcal{O}_C(P_{N^2-1}+P_{N^2})$ moves in the family $A\otimes \mathcal{O}_C(P+Q)$ with $P,Q$ moving in $C$. None of this line bundles but $A\otimes \mathcal{O}_C(P_{i}+P_{j})$ come from $\mathbb{P}(F_2)$, otherwise, there would be two different maps $S\to \mathbb{P}^1$ inducing the same on $C$. We get that $\mathrm{Im}( f_2)\cap W^1_{g(S,L)}(C)$ is zero dimensional at $A\otimes \mathcal O_C(P_{N^2-1}+P_{N^2})$. This implies that the localization of the fiber of $\pi$ above $C$ at is of dimension zero, hence we get the claim. The Theorem follows. 
  
  \end{proof}
   Notice that  $\mathrm{deg}({{\mathrm{Im}(f_2)}\to |L|})\ge \frac{N^2(N^2-1)}{2}$, so that  the general curve in $|L|$ carries at least $\frac{N^2(N^2-1)}{2}$ minimal pencils. We conjecture that there are exactly this number of minimal pencils if $S$ has rank $1$ Néron-Severi.
We may finally prove (in the hypothesis of Theorem \ref{t1}):
\begin{corollary}
  $\mathrm{dim}\mathcal{W}^1_{g(S,L)}|L|\ge g-2$ and $\mathrm{gon(C)} \ge g(S,L)$. Moreover, the general curve in $|L|$ carries at least $\frac{N^2(N^2-1)}{2}$ minimal pencils.
  \label{whatuneed}
\end{corollary}
We have finally established Theorem \ref{t1}
\begin{corollary}
  The linear growth condition (see \cite{Ap}, Theorem 2) holds for the general curve in $|L|$ where $L$ is a linear system on a simple abelian surface satisfying the hypotheses of Theorem \ref{t1}.
\end{corollary}
\begin{proof}
  On the one hand by Theorem \ref{t3} we have $\mathrm{dim}\mathcal{W}^1_{g(S,L)}|L| \le g-2$ and $\mathrm{dim}\mathcal{W}^1_{g(S,L)-1}|L| \le g-3$. On the other hand, we just proved (corollary \ref{whatuneed}) that each smooth curve in $|L|$ carries at least one $g^1_{g(S,L)}$ so that $\mathrm{dim}\mathcal{W}^1_{g(S,L)}|L| \ge g-2$. It follows that the gonality of the general curve in $|L|$ is $g(S,L)$. The linear growth condition follows again by Theorem \ref{t3}.

  The claim about the loci of curves of lower gonality follows with an analogous argument.
  \end{proof}

Let us also notice that for the general curve $C \in L$ the locus $W^1_{g(S,L)}(C)$ is non-reduced.
\begin{corollary}
  For general $C$ the Brill--Noether locus $W^1_{g(S,L)}(C)$ is finite and non-reduced at all its points. Moreover, the general curve carries at least $\frac{N^2(N^2-1)}{2}$ minimal pencils. 
  \end{corollary}
   
\begin{proof}
   Let $(C,A)$ be a general point in the image of $\mathbb{P}(F_2) \to \mathcal{W}^1_{g(S,L)}|L|$. We have  $H^0(N)=H^0(N_{|_{C}})$ and $H^0(N^{\otimes 2})=H^0(N^{ \otimes 2}_{|_{C}})$ (Since $L>2N$ we have $h^1(L \otimes N^{\otimes -2})=0$). So $H^0(A^{\otimes 2})$ are just sections of $H^0(N^{\otimes 2})$ tangents to $C$ at $P_1,\dots,P_{N^2-2}$, let us call this vector space $V$. We have $h^0(A^{\otimes 2})= \mathrm{dim} V \ge 4$, indeed we are imposing $4h^0(N)-4$ linear conditions on $H^0(N^{\otimes 2})$ (which has dimension $4h^0(N)$). By Brill--Noether theory the tangent space of $W^1_{g(S,L)}$  at $A$ is given by the image of $H^0(A^{\otimes 2}) \to H^0(\omega_C)^\vee$. Hence of dimension $h^0(A^{\otimes 2})-3=1$. By the previous corollary $\mathrm{dim}(W^1_{g(S,L)}(C))=0$. Hence  $W^1_{g(S,L)}(C)$ is non-reduced at $A$.
  A different proof may be given by using Proposition \ref{313}. This second proof also gives more evidence for Donagi--Morrison conjecture on abelian surfaces. Indeed, it says that for the general curve the Brill--Noether locus is non-reduced at all its points (hence all of them may, in principle, be points coming from $\mathbb{P}(F_2)$).
  \end{proof}

    \section{On certain special loci in $\mathrm{Hilb}(S)$}
      The goal of this section is to prove the following:
      \begin{proposition}
  Write
    \[
  \mathcal{P}_{[n],l}= \bigcup_{i=0}^4 \mathcal{Z}^i_{[n],l},
  \]
  where
  \begin{eqnarray}
    \mathcal{Z}^i_{[n],l}=\{ E \in \mathcal{P}_{[n],l} \textrm{ so that } h^1(E)=i; \} \nonumber
    \end{eqnarray}
    then for general $E \in  \mathcal{Z}^i_{[n],l}$ the image of the projection $q:(\mathcal{Z}^i_{[n],l})_{(N,\eta)} \to \mathrm{Hilb}^l(S)$ is at most $(2l-2i)$-dimensional at $q(E)$.
  \label{diff}
  \end{proposition}
      First, we sketch the proof so that we introduce some notation and the main tool we need.
      \begin{proof}[sketch]
        Let \[
      A_{l,(N,\eta),i}= q((\mathcal{Z}^i_{[n],l})_{(N,\eta)}).
      \]
      Let us consider the Grassmannian bundle $B_{l,(N,\eta),i} \to A_{l,(N,\eta),i}$, where
      \[
      (B_{l,(N,\eta),i})_\xi=\mathrm{Gr}(2,H^0(N \otimes \mathcal{I}_\xi)).
      \]
      Now we consider
      \[
      \varphi: B_{l,(N,\eta),i} \dashrightarrow \mathrm{Hilb}^l(S),
      \]
      where $\varphi(V)=Z(s_1) \cap Z(s_2)$ is the intersection scheme of two sections generating $V \subset H^0(N\otimes \mathcal{I}_\xi)$. Since for any Lazarsfeld--Mukai bundle $|N \otimes \mathcal{I}_\xi|$ has zero dimensional base locus the above map is generically well defined. Moreover, the image of $\varphi$ is contained in $\mathrm{Gr}(2,H^0(N)) \dashrightarrow \mathrm{Hilb}^{N^2}(S)$. Now if $\xi$ is reduced  for $E$ general in its irreducible component it is clear that the above map has finite fibers. From an easy computation one gets the desired inequality.
   In order to complete the proof one must analyze carefully the components of $ B_{l,(N,\eta),i}$ when $\varphi$ has not generically finite fibers. We will prove in this section that the locus where the fibers are of dimension $l$ is of codimension $\ge l$ in the image of $\mathrm{Gr}(2,H^0(N)) \to \mathrm{Hilb}^{N^2}(S)$.
      \end{proof}

      We need a number of preliminary lemmas. It is worth mentioning here that the above proof with a minimal extra argument already works for $i=1$ which is the needed value to bound the dimension of $\mathcal{W}^1_{g(S,L)-1}|L|$, which is sufficient in order to establish that the gonality is at least $g(S,L)$. Moreover, the arguments in the appendix may be simplified a lot with some positivity assumption on $N$ (i.e. $N^2 >>0$), see the remark below lemma \ref{mull}. From now on $\zeta \in \mathrm{Im}(\varphi)$ and $(\xi,V)$ is a general element of $\varphi^{-1}(\zeta)$. Moreover, we will fix one for all an irreducible component of $\mathcal{G}\subset \mathcal{G}_{[n],l}$ dominating a component $W \subset \mathcal{W}^1_d|L|$ via the construction described before  such that $W \to |L|$ is dominant (see the beginning of section 2.2). Moreover (abusing notation) $ A_{l,(N,\eta),i}$ will denote just the image of $(\mathcal{Z}^i_{[n],l})_{(N,\eta)}$ intersected with the chosen component in $\mathrm{Hilb}^lS$. As before $gr:B_{l,(N,\eta),i} \to A_{l,(N,\eta),i}$ will denote the relative Grassmannian of planes.  We are interested in the dimension of $\varphi^{-1}(\zeta) \subset B_{l,(N,\eta),i}$, one obvious condition is that $\xi \subset \zeta$ and $\mathcal{I}_\zeta \subset \mathcal{I}_\xi$. The first lemma we need is a refinement of this condition.
      \begin{lemma}
        For general $(\xi,V) \in  B_{l,(N,\eta),i}$ for any $P \in S$ there exists $u \in \mathcal{O}_{S,P} -m_P^2$ such that $u \mathcal{I}_\xi \subset \mathcal{I}_\zeta$. Moreover, $u$ can be chosen locally to be the equation of a smooth curve $C$ such that there exists $A \in W\cap W^1_d(C)$ whose associated Lazarsfeld--Mukai bundle may be realized as an element of $\mathrm{Ext}^1(N \otimes \mathcal{I}_\xi, M \otimes \eta)$.
        \label{abv}
        \end{lemma}
      \begin{proof}
        We use the hypothesis that they come from a Lazarsfeld--Mukai bundle associated to a linear series on a smooth curve. We may lift $s_1,s_2$ (generators of $V$) to two sections of $E$ such that $u=s_1 \wedge s_2$ gives the equation of a smooth curve. Now it is routine to verify that for any $\xi$ we have $ u\mathcal{I}_\xi \subset \mathcal{I}_\zeta \subset \mathcal{I}_\xi$. Let us check this \'etale locally. The local expressions for $s_1,s_2$ are
      \begin{eqnarray}
        s_1=(f_1,g_1),s_2=(f_2,g_2) \in \mathbb{C}[[x,y]] \times \mathbb{C}[[x,y]], \nonumber
        \end{eqnarray}
      so we get $s_1 \wedge s_2=f_1g_2-g_1f_2$. The local expression for $E \to N  \otimes \mathcal{I}_\xi$ is of the form
      \begin{eqnarray}
        (1,0) \to p, \nonumber \\
        (0,1) \to q;  \nonumber
      \end{eqnarray}
      with $p,q \in \mathbb{C}[[x,y]]$. So we get $\zeta=(pf_1+qg_1,pf_2+qg_2)$. Now we obtain $pu=g_2(pf_1+qg_1)-g_1(pf_2+qg_2)$ and $qu=-f_2(pf_1+qg_1)+f_1(pf_2-qg_2)$. The key point is that $u$ defines a smooth curve so $u \in  \mathcal{O}_{S,P} -m_P^2$. 
      \end{proof}
      \begin{rmk}
        It is clear that the source of the positive dimension fibers of $\varphi$ is the non-reduceness of $\xi$. Actually, the more it is true: if the fibers are positive-dimensional it means that the general Lazarsfeld--Mukai bundle lying above $\xi$ is associated to a smooth curve passing through at least one of the points where the localization of $\xi$ is non-reduced.
        \end{rmk}
      Let us deduce some corollaries.
      \begin{corollary}
          $(\mathcal{I}_\xi)_P$ is principal in $\mathcal{O}_{S,P}/(\mathcal{I}_\zeta)_P$.
      \end{corollary}
      \begin{proof}
        The fact that $(\mathcal{I}_\xi)_P$ is a principal ideal in $\mathcal{O}_{S,P}/(\mathcal{I}_\zeta)_P$ follows from the smoothness of $u=f_1g_2-g_1f_2$. Indeed, we get that at least one among $f_1,f_2,g_1,g_2$ is invertible at $P$; without loss of generality $f_1 \in \mathcal{O}_{S,P}^*$. Then, $p=f_1^{-1}((pf_1+qg_1)-qg_1) \in (q,\mathcal{I}_{\zeta_P})$.
      \end{proof}
      \begin{corollary} Write $(\zeta)_P=(p,q)$ and $(\xi)_P=(a,b)$ in the \'etale local ring at $P$ (they both are always local complete intersection). Then (up to interchange $p$ and $q$)  $\mathrm{mul}(a)-1 \le\mathrm{mul}(p) \le \mathrm{mul}(a)$
        \label{abpq}
      \end{corollary}
      \begin{proof}
     Immediate from the previous lemma.
        \end{proof}
 
    \begin{lemma}
     The locus in the Hilbert scheme of length $6$ subschemes of of $\mathrm{Spec} \ \mathbb{C}[[x,y]]$ whose ideal is a complete intersection $(p,q)$ with $\mathrm{mul}(p)=2,\mathrm{mul}(q)=3$ is $4-$ dimensional.
      \label{H}
    \end{lemma}
    \begin{proof}
      Let $V_i$ be the vector space of homogeneous polynomials of degree $i$.
      Let us consider the projective bundle $W \to V_2-\{0\}$ where $W_g=\mathbb{P}(V_3/\langle xg,yg \rangle)$ equipped with the tautological line bundle $F \to W$ where $F_{g,w}=V_3/\langle xg,yg,w \rangle$. We have the rational morphism
      \[
      \begin{tikzcd}
        F \arrow{r}{\varphi} & \mathrm{Hilb}^6 (\mathrm{Spec} \ \mathbb{C}[[x,y]]),
        \end{tikzcd}
      \]
      given by $(g,w,f) \to \mathbb{C}[[x,y]]/(g+f,w)$. It is easy to see that the morphism is well defined (outside the locus of non transverse intersection) and dominates the special locus in the Hilbert scheme we are interested in. The fibers of $\varphi$ are $1-$ dimensional since $\varphi(g,w,f)=\varphi(\lambda g, w, \lambda f)$ for any $\lambda \in \mathbb{C}^*$. The lemma follows.
      \end{proof}
    \begin{lemma}
      If an ideal $\mathcal{I} \subset \mathbb{C}[[x,y]]$ is saturated at level $n$ (i.e. $(x,y)^n \subset \mathcal{I}$) and there exists $f \in \mathcal{I}$ with $\mathrm{mul}_0(f)=2$, then $\mathcal{I}$ contains either a nodal curve or a cuspidal curve of order $n$ (we mean $x^2=y^n$ for $n \ge 3$ up to change of coordinates).
      Moreover, the principal ideals of a cusp $\mathcal{I} \subset \mathbb{C}[[x,y]]/(x^2-y^n)$ are of the following shape:
       \begin{eqnarray}
        (xy^l+\lambda_1 y^{l+1}+\dots+\lambda_{n-1}y^{l+n-1}), \quad \quad \quad \lambda_i \in \mathbb{C} ;\nonumber \\
        (y^{l+1}).   \ \ \ \ \quad \quad \quad \quad \quad \quad \quad \quad \quad \quad \quad \quad \quad \quad \quad \quad \quad  \nonumber
      \end{eqnarray}
   
     \label{cusp} 
    \end{lemma}
    \begin{proof}
      If $f$ is nodal then we are done.
      Otherwise, since $\mathcal{I}$ is saturated at level $n$ we may suppose $f=x^2w_1+w_2y^k+xr(y)=xw(xw+r(y)w^{-1})+y^kw_2$, with $w^2=w_1,w_2$ invertible (an invertible element always admits a $n$-th root in a power series ring). Now let $x'=wx+\frac{r(y)}{2w}$, we obtain $f=x'^2-\frac{r(y)^2}{4w_1}+y^kw_2=x'^2+y^sw_3$, where $s=\mathrm{min}\{k,2\mathrm{mul}_0(r)\}$, $w_3$ is invertible. Up to a change of coordinates, we got the desired result. The facts concerning the principal ideals are easy computations.
      
      \end{proof}

        The following corollary will be the key in studying the fibers of $\varphi$.
      \begin{corollary}
          Let us suppose \'etale locally $\mathcal{I}_{\zeta_P}=(x^2-y^n,xy^l+\lambda_1 y^{l+1}+\dots+\lambda_{n-1}y^{l+n-1}) $. Then the following hold:
          \begin{itemize} 
              \item The locus of principal subschemes in $\mathrm{Hilb}^8(\mathrm{Spec} \, \mathbb{C}[[x,y]]/\mathcal{I}_{\zeta_P})$ is at most $4-$dimensional.
              \item The locus of principal subschemes in $\mathrm{Hilb}^j(\mathrm{Spec} \, \mathbb{C}[[x,y]]/\mathcal{I}_{\zeta_P})$ is at most $3-$dimensional for $j=6,7$.
              \item The locus of principal subschemes in $\mathrm{Hilb}^6(\mathrm{Spec} \, \mathbb{C}[[x,y]]/\mathcal{I}_{\zeta_P})$ whose associated ideal is of the form $(xy^k +\delta_1 y^{k+1}+\dots+\delta_{n-1}y^{k+n-1})$ with $k \ge 1$ is at most $2-$dimensional.
              \item The locus of principal subschemes in $\mathrm{Hilb}^j(\mathrm{Spec} \, \mathbb{C}[[x,y]]/\mathcal{I}_{\zeta_P})$ is at most $2-$dimensional for $j \le 5$.
          \end{itemize}
          \label{ganga}
              \end{corollary}
          \begin{proof}
         We give a sketch.
              We work modulo $\mathcal{I}_{\zeta_P}$. Let us take $\xi \subset \zeta_P$. Since $\mathrm{dim}_\mathbb{C}\mathbb{C}[[x,y]]/\mathcal{I}_{\xi}=8$ we may suppose $y^8 \in \mathcal{I}_\xi$. Let us suppose $\mathcal{I}_\xi=(x+\delta_1y+\dots+\delta_{n-1}y^{n-1})$, we get $\delta_1,\delta_2,\delta_3\in\{0;-1\}$ (otherwise $y^6\in \mathcal{I}_\xi$ and  $\mathrm{dim}_\mathbb{C}\mathbb{C}[[x,y]]/\mathcal{I}_{\xi}\le 6$). For istance if $\delta_1=\delta_2=\delta_3=0$, we get $\mathcal{I}_\xi=(x+\delta_1y+\dots+\delta_{n-1}y^{n-1})=(x+\delta_4y^4+\delta_5y^5+\delta_6y^6+\delta_7y^7)$ depends on (at most) $4-$parameters (the other cases being similar). Now suppose  $\mathcal{I}_\xi=(xy+\delta_1y^2+\dots+\delta_{n-1}y^{n})$. Since $\{1,x,y\}$ are linearly indipendent in $\mathbb{C}[[x,y]]/\mathcal{I}_{\xi}$  and $\mathrm{dim}_\mathbb{C}\mathbb{C}[[x,y]]/\mathcal{I}_{\xi}=8$ we deduce that $y^7 \in \mathcal{I}_\xi$. Moreover $\delta_1,\delta_2 \in \{0;-1\}$ (otherwise $y^5 \in \mathcal{I}_\xi$ and  , we get that up to the discrete data of $\delta_1,\delta_2$ the ideal  $\mathcal{I}_\xi$ depends only on $\delta_3,\delta_4,\delta_5$ (at most $3-$parameters). The other cases and points can be proven with a similar reasoning.
          \end{proof}

   Before finally takling the proof of \ref{diff} we need some more lemmas to bound the dimension of some subschemes of $\mathrm{Gr}(2,H^0(N))$. We report also a proof for a lack of reference.
    \begin{lemma}
      Let $S$ be a simple abelian surface and $N \in \mathrm{Pic}(S)$ an ample line bundle. Then the following hold:
      \begin{enumerate}
           \item For a general point $P \in  S$ the codimension of $H^0(N \otimes \mathcal{I}_P^2) \subset H^0(N)$ is  $3$ and the codimension of $H^0(N \otimes \mathcal{I}_P^3)\subset H^0(N)$ is at least $4$;
      \item  For a general point $P \in T \subset  S$ the codimension of $H^0(N \otimes \mathcal{I}_P^2) \subset H^0(N)$ is  at least $2$ and the codimension of $H^0(N \otimes \mathcal{I}_P^3)\subset H^0(N)$ is at least $3$, where $T \subset S$ is a one dimensional subscheme;
   \item For any $P \in S$ the codimension of $H^0( N \otimes \mathcal{I}_P^2)$ is at least $1$;
        \item For $(P,Q)$ general in a codimension $i=0,1$ subscheme $T \subset S^2$ the codimension of \\ $H^0(N \otimes \mathcal{I}_{P,Q}^2) \subset H^0(N)$ is at least $6-i$.
      \end{enumerate}
      The codimension of a vector space of dimension $0$ or $1$ is considered to be $\infty$.
      \label{mull}
    \end{lemma}
    \begin{proof}
      We prove the first point the other ones are similar.
      For the facts concerning the codimension of $H^0(N \otimes \mathcal{I}_P^2)$ just consider the morphism induced by $N$. As long as $h^0(N) \ge 3$ the image is a surface, hence the map is locally \'etale out of the base points (finite) and the ramification locus (one dimensional subscheme $R \subset S$). It is easy to verify that for the general point  $P \in S-R$ we have  $H^0(N \otimes \mathcal{I}_P^2)$ of codimension $3$ in $H^0(N)$. Indeed, take a smooth curve $C\subset S$ at $P$ (it exists for general $P$), then by local \'etalness we have $H^0(N \otimes (\mathcal{I}_P^2,u))\subset H^0(N)$ of codimension $2$ (the morphism separates tangent directions). Now, since the equation $s$ of $C$ lies in $H^0(N \otimes (\mathcal{I}_P^2,u))-H^0( N \otimes \mathcal{I}_P^2)$ we obtain the desired result. Now we deal with $H^0(N \otimes \mathcal{I}_P^3)$ for $P$ general in $S$. Suppose by contradiction there exist $P_1, \dots , P_{\lfloor \frac{N^2-1}{6} \rfloor}$ general distinct points such that $H^0(N \otimes \mathcal{I}_P^3) \subset H^0(N)$ is of codimension $ \le 3$. We get $H^0(N \otimes \mathcal{I}_{P_1,\dots,P_{\lfloor\frac{N^2-1}{6} \rfloor}})$ is of dimension $\ge 1$. Now for $N^2 \ge 14$ and $N^2 \neq 18$ if we take a curve in $C_1 \in |H^0(N \otimes \mathcal{I}_{P_1,\dots,P_{\lfloor\frac{N^2-1}{6} \rfloor}})|$  and a curve $C_2 \in |H^0(N \otimes \eta \otimes \mathcal{I}_{P_1, \dots , P_{\lfloor\frac{N^2-1}{6} \rfloor}})|$ for $\eta$ non trivial in $\mathrm{Pic}^0(S)$ we get $C_1\cdot C_2 >N^2$ contradiction. If $N^2=12,18$ it is easy to see that at a general point $P \in S$ there is a curve of multiplicity $2$ (either using $N=N_1+N_2$, or if $N$ is primitive by taking an \'etale cover $(S,N) \to (\tilde{S},\tilde{N})$ with $\tilde{N}$ of type $(1,3)$ and considering the preimage of a singular curve in $|H^0(\tilde{N})|$). Hence $H^0(N \otimes \mathcal{I}_P^3) \subset H^0(N \otimes \mathcal{I}_P^2)$ is of codimension at least $1$ and the result follows by what we proved above. For $N^2=10$ we fix $C_0 \in |H^0(N \otimes  \eta \otimes \mathcal{I}_{P_1}^3)|$ for $\eta \in \mathrm{Pic}^0(S)- \mathcal{O}_S$. Then we have  a non trivial map $\mathbb{P}^1 \simeq |H^0(N \otimes \mathcal{I}_{P_1}^3)|\to S$ given by  $C \to C_0\cdot C-9P$, contradiction. For $N^2 <10$ there is nothing to prove by numerical reasons.
      \end{proof}

    We finally take care of the proof of \ref{diff}.
    \begin{proof}[Proof of \ref{diff}]
      We will consider just the case $i=4$ which is the most complicated and leave the others ($i<4$) to the reader. Let us observe that the Hilbert scheme of points supported at one point on a smooth curve is discrete, hence if $(\mathcal{I}_\zeta)_P$ contains a polynomial defining a smooth curve then the fibers of $\varphi$ are finite and there is nothing to prove. Hence, we may suppose there is a point $P$ such that $\mathcal{I}_{\zeta_P}=(a,b) \subset\mathcal{I}_P^2$ and $\mathrm{length}(\xi_P) \ge 2$.
      We split the proof in several cases according to the most special point in $\mathrm{Supp}(\xi)$.
      \begin{enumerate}
      \item First we suppose that for general $\xi \in A_{l,(N,\eta),4}$ there exists $P \in S$ such that $\mathrm{length}(\xi)_P \ge 8$. We get that $P(\xi)$ is the unique non-reduced point and that $P$ is general in $S$ otherwise $A_{l,(N,\eta),4}$ is contained in a codimension $8$ subscheme and we are done. So, the only contribution to positive-dimensional fibers comes from the variation of $\xi_P$ as a subscheme of $\zeta_P$. Now write $\xi=(p,q),\zeta=(a,b)$. Observe that by \ref{abpq} we must have $\mathrm{mul}_{P(\xi)}(a) \le 3$ (up to interchange $a$ with $b$). We split the proof in two cases according to the multiplicities  of $a,b$. The possibilities we have to study are given by the fact that $\mathrm{length}(\xi_P)=8$ (so $\mathrm{mul}(p)\mathrm{mul}(q) \le 8$) together with lemma \ref{abpq}.
        \begin{enumerate}
        \item Case $\mathrm{mul}(a)=3,\mathrm{mul}(b) \ge 3$. In this case the image of $B_{l,(N,\eta),4} \to \mathrm{Hilb}^{N^2}(S)$ is generically contained in the image of
          \[
          \bigcup_{P \in U} \mathrm{Gr}(2,H^0(N\otimes \mathcal{I}_P^3)) \to \mathrm{Hilb}^{N^2}(S),
          \]
          where $U$ is the open subset where $H^0(N \otimes \mathcal{I}_P^3)$ is of minimal dimension. $H^0(N \otimes \mathcal{I}_P^3)\subset H^0(N)$ is of codimension at least $4$ for general $P$ by lemma \ref{mull}. So $\bigcup_{P \in U} \mathrm{Gr}(2,H^0(N\otimes \mathcal{I}_P^3))$ is of dimension at most $2(h^0(N)-6)+2=2h^0(N)-10$. Now the fibers $B_{l,(N,\eta),4} \to \mathrm{Hilb}^{N^2}(S)$ are of dimension at most $6$; indeed we must have $\mathrm{mul}(p)=2$ and the quadratic part of $p$ must divide the cubic part of something of the form $\lambda a +\delta b$ with $\gamma,\delta \in \mathbb{C}$. This gives a locally codimension $1$ condition on the $7$-dimensional Hilbert scheme of points of length $8$ supported at a given point. This gives us $\mathrm{dim}(B_{l,(N,\eta),4}) \le 2h^0(N)-4$ which is what we want.
        \item Case $\mathrm{mul}(a)=2$.  In this case the image of $B_{l,(N,\eta),4} \to \mathrm{Hilb}^{N^2}(S)$ is generically contained in the image of
          \[
          \bigcup_{P \in U} \mathrm{Gr}(2,H^0(N\otimes \mathcal{I}_P^2)) \to \mathrm{Hilb}^{N^2}(S),
          \]
          where $U$ is the open subset where $H^0(N \otimes \mathcal{I}_P^2)$ is of minimal dimension. Moreover, by lemma \ref{cusp} we may suppose there is a cusps in $\mathcal{I}_{\zeta_P}$. Hence, by corollary \ref{ganga} the fibers are of dimension at most $4$. The Proposition follows with a computation analogous to the previous one.
         
   \end{enumerate}
      \item  Now suppose that for general $\xi \in A_{l,(N,\eta),4}$ there exists $P \in S$ such that $\mathrm{length}(\xi)_P \ge 7$. If $P(\xi)$ is general in $S$ then one can basically argue as in point $1$ (the dimension of the fibers due to the point of length $7$ drops off by one, but there may be another non-reduced point in the support of $\xi$ making them jump by one again; we leave the details to the reader). If $P(\xi) \in C$ a curve in $S$ then one has to be a bit more careful, but the ideas are again the same as in (1). In this case we may suppose the general $\xi$ is reduced outside $P(\xi)$ (otherwise we would be already in a codimension $8$ subscheme and we would be finished). The key observation to reduce to (1) is observing that in this case the fibers drops always dimension by one because the length is one less and the dimension of the subschemes of $\mathrm{Hilb}^{N^2}(S)$ we need to consider increases at most by $1$ by lemma \ref{mull} and considerations on the cusps similar to the previous point. We leave the details again to the reader.
 
\item  Now suppose that for general $\xi \in A_{l,(N,\eta),4}$ there exists $P \in S$ such that $\mathrm{length}(\xi)_P \ge 6$. If $P(\xi)$ is general in $S$ then one can basically argue as in point $1$ (the dimension of the fibers due to the point of length $6$  drops off by two, but there may be two more non-reduced points in the support of $\xi$ making them jump by two again, we leave the details to the reader). If $P(\xi) \in C$ a curve in $S$ then one can argue basically as in (2). If $P(\xi)=P$ a fixed point in $S$ varying $\xi \in A_{l,(N,\eta),4}$ one has to be a bit more careful, but the ideas are again the same as in (1). Let us provide some more details in this last situation. In this case we may suppose the general $\xi$ is reduced outside $P(\xi)$ (otherwise we would be already in a codimension $8$ subscheme and we would be finished).

  \begin{enumerate}
   \item  Case $\mathrm{mul}(a)=3,\mathrm{mul}(b) \ge 3$. This forces $(p,q)$ to be of multiplicities $2,3$ respectively and of complete intersection. A quick parameter count yields that the subschemes of this shape describe a $4-$ dimensional subscheme of $\mathrm{Hilb}^6(\mathrm{Spec} \ \mathbb{C}[[u,y]])$ (see lemma \ref{H}). Hence $A_{l,(N,\eta),4} \subset \mathrm{Hilb}^l(S)$ is of codimension $8$, via direct computation.
   \item Case $\mathrm{mul}(a)=2$. In the following $u=u(\xi)$. The coordinates $u,y$ will be local coordinate at $P$ changing accordingly to $\xi$ ($u$ will denote the tangent direction of the cusps in $\zeta_P)$. If the ideal of $\xi_P \subset \zeta_P$ is cut by an equation of the form $(uy^k+\delta_1 y^{k+1}+\dots+\delta_{n-1}y^{n+k-1})$ with $k \ge 1$ in the notation of corollary \ref{ganga}, then, by the aforementioned corollary, the fibers of $\varphi$ are of dimension at most $2$ and a computation similar to the first point proves the proposition. If the ideal of $\xi_P \subset \zeta_P$ is cutted by an equation of the form $(u+\delta_1 y+\dots+\delta_{n-1}y^{n-1})$ then either the equation given by lemma \ref{abv} has $u$ as a linear term or $\zeta_P$ contains a nodal curve (hence the fibers of $\varphi$ are $1-$ dimensional). So we may suppose the equation given by lemma \ref{abv} lies in $(\mathcal{I}_p^2,u)$ and that the only curves of multiplicity $2$ in $\zeta_p$ have quadratic term $u^2$. Now the locus $A_{u_0}=\{ \xi \in \mathrm{Hilb}^6(S)|$ $\mathcal I_{\xi_P}\xi_P$ contains a smooth curve with direction $u_0\}$ is of codimension $8$. So we may suppose $A_{u_0} \subset  A_{l,(N,\eta),4}$ of codimension $\ge 1$, which implies in particular that at $P$ we have cusps with any tangent direction. We get $H^0( N \otimes (\mathcal{I}_P^3,u^2)) \subset H^0(N)$ of codimension $\ge 3$. So, the image of $\varphi$ is generically contained in
     \[
 \bigcup_{u \in \mathbb{P}({T_P^\vee S})} \mathrm{Gr}(2,H^0(N\otimes (\mathcal{I}_P^3,u^2))) \to \mathrm{Hilb}^{N^2}(S),
 \]
 a computation similar to that of point $1$ completes the proof.
  \end{enumerate}
 \item  Now suppose that for general $\xi \in A_{l,(N,\eta),4}$ there exists a unique point $P \in S$ such that $\mathrm{length}(\xi)_P \ge 2$ and $\mathcal{I}_{\zeta_{P_i}}\subset \mathcal{I}_{P}^2$. The arguments we explained in the previous points generalize to this case easily.
 \item In order to end the proof we deal with the case $\xi \in A_{l,(N,\eta),4}$ has (at least) two non-reduced points giving rise to positive dimensional fibers (this implies that $\mathrm{length}(\xi_{P_i}) \ge 2$ and $\mathcal{I}_{\zeta_{P_i}}\subset \mathcal{I}_{P_i}^2$). Then we have a natural map $(P_1+\dots+P_n)(\xi) \in \mathrm{Sym}^n(S)$. Let us observe that we may suppose \[
   \sum_{i=1}^n(\mathrm{length}(\xi_{P_i})-1)+ \mathrm{codim}((P_1+\dots+P_n)(A_{l,(N,\eta),4}))\subset \mathrm{Sym}^n(S)) \le 8
   \] otherwise $A_{l,(N,\eta),4}$ would be contained in a codimension $8$ closed subscheme of $\mathrm{Hilb}^l(S)$ which is what we want. Hence the fibers of $\varphi$ are of dimension at most $7$. Moreover we may pullback the locus $(P_1+\dots+P_n)(A_{l,(N,\eta),4})$ to $A \subset S^n$. We fix a component of $A$ (all the components are isomorphic via the action of the symmetric group).
   \begin{enumerate}
   \item If $(P_i,P_j) \in S^2$ is of codimension $k$ for some $(i,j)$ (meaning that the image $T$ of $(pr_i,pr_j):A \to S^2$ of codimension $k$) then the image of $\varphi$ is generically contained in
     \[
     \bigcup_{(P_1,P_2) \in U} \mathrm{Gr}(2, H^0(N \otimes (\mathcal{I}_{P_1,P_2}^2))) \to \mathrm{Hilb}^{N^2}(S)
     \]
     which is of dimension at most $2(h^0(N)-8+k)+4-k=2h^0(N)-12+k$ ($U\subset T$ is a generic open set where $H^0(N \otimes (\mathcal{I}_{P_1.P_2}^2))$ has maximal codimension $\ge 6- k$ by lemma \ref{mull}). The fibers of $\varphi$ are of dimension at most $7-k$ (since $\sum_{i=1}^n\mathrm{length}(\xi_{P_i})+ \mathrm{codim}((P_1+\dots+P_n)(A_{l,(N,\eta),4}))\subset \mathrm{Sym}^n(S))\le 8$). The Proposition follows.
     \item If $(P_i,P_j) \in S^2$ is of codimension $k>1$ for some $(i,j)$ (meaning that the image $T$ of $(pr_i,pr_j):A \to S^2$ of codimension $k$) then the image of $\varphi$ is generically contained in
     \[
     \psi:\bigcup_{(P_1,P_2) \in U} \mathrm{Gr}(2, H^0(N \otimes (\mathcal{I}_{P_1,P_2}^2))) \to \mathrm{Hilb}^{N^2}(S)
     \]
     which is of dimension at most $2(h^0(N)-8+k+1)+4-k=2h^0(N)-10+k$ ($U\subset T$ is a generic open set where $H^0(N \otimes (\mathcal{I}_{P_1.P_2}^2))$ has maximal codimension $\ge 5- k$ by lemma \ref{mull}). Now, if the fibers are of dimension $<7-k$ we are done (with a computation analogous to the previous point). Otherwise, we must have $\mathrm{length}(\xi_{P_1})=2+i,\mathrm{length}(\xi_{P_2})=7-k-i$ with $i=0,1$. Moreover, on the one hand $\varphi(\xi,V)$ must contain only cusps at the point $P_2$ of length $\ge 3$ (otherwise the dimension of the fibers of $\varphi$ would drop since nodal curves have $1$-dimensional Hilbert scheme of any degree). On the other hand,  the condition that $\varphi(\xi,V)$ contains only cusps with a given tangent direction imposes a codimension one condition on subschemes of length $\ge 3$ of $\mathrm{Spec} \ \mathbb{C}[[x,y]]$ (so the dimension of the fibers would drop again). The Proposition follows.
     
     \end{enumerate}
          \end{enumerate}
     
        \end{proof}
    
\end{document}